\documentclass[twoside,reqno,11pt]{amsart}
\usepackage{upref,amsxtra,amssymb,amscd}
\usepackage{varioref}
\usepackage{verbatim}
\usepackage{epsfig}
\usepackage{psfrag}
\usepackage{graphicx}
\usepackage{epic,eepic,eucal}
\usepackage{enumerate}

\def\ci{     \circ} 
 
\def\A{          \mathcal A}
\def\mS{          \mathcal S}

 \def\a{         \alpha}

\def\de{ \delta  }
\def\si{ \sigma   }
\newcommand{\AAA}{{\mathbb A}}
\newcommand{\RR}{{\mathbb R}} 
\newcommand{\LL}{{\mathbb S}}
\newcommand{\NN}{{\mathbb N}}
\newcommand{\TT}{{\mathbb T}}
\newcommand{\ZZ}{{\mathbb Z}}
\newcommand{\CC}{{\mathbb C}}
\newcommand{\QQ}{{\mathbb Q}}
\newcommand{\DD}{{\mathbb D}}

\def\carre{ \hfill $\Box$    }

\DeclareMathOperator{\diff}{Diff}

\newtheorem{theo}{\sc Theorem}[section]
\newtheorem{prop}[theo]{\sc Proposition}
\newtheorem{lemm}[theo]{\sc Lemma}
\newtheorem{coro}[theo]{\sc Corollary}

\theoremstyle{definition}

\newtheorem{defi}[theo]{\sc Definition}

\theoremstyle{remark}

\newtheorem{rema}[theo]{\sc Remark}

\numberwithin{equation}{section}

\begin{document}

\newcommand{\la}{\lambda}
\newcommand{\nbd}{neighborhood}
\newcommand{\e}{\varepsilon}
\newcommand{\ph}{\varphi}
\newcommand{\Ph}{\Phi}
\newcommand{\ti}{\tilde}
\newcommand{\De}{\Delta}
\newcommand{\Diffr}{\text{Diff}^\omega_r(\Bbb T^2)}
\newcommand{\Diffk}{\text{Diff}^k(\Bbb T^2)}
\newcommand{\Diffi}{\text{Diff}^{\infty}(\Bbb T^2)}
\newcommand{\be}{\beta}
\newcommand{\ga}{\gamma}
\newcommand{\I}{\text{Im\,}}

\newcommand{\pdvr}[2]
{\dfrac{\partial^{#2} #1}{\partial \theta^{#2_1} \partial r^{#2_2}}}
%
\newcommand{\pdvx}[2]
{\frac{\partial^{#2} #1}{\partial x^{#2_1} \partial y^{#2_2}}}
\newcommand{\pdvrs}[2]
{\partial^{#2} #1 /\partial \theta^{#2_1} \partial r^{#2_2}}
%
\newcommand{\pdvxs}[2]{
\partial^{#2} #1/\partial x^{#2_1} \partial y^{#2_2}}

\title[Weak mixing disc and annulus diffeomorphisms]{Weak mixing disc and annulus diffeomorphisms with 
arbitrary Liouville rotation number on the boundary.}

\author{B.~Fayad, M.~Saprykina}
\address{Bassam  Fayad, Universit\'e Paris 13, CNRS 7539, Villetaneuese 93430  }
\email{fayadb@math.univ-paris13.fr}
\address{Maria Saprykina, Department of Mathematics, University of Toronto,
Toronto, ON, Canada M5S 3G3}
\email{masha@math.kth.se}

\maketitle 


\begin{abstract}

Let $M$ be an $m$-dimensional differentiable manifold with a nontrivial circle action ${\mathcal S}= 
{\lbrace S_t \rbrace }_{t \in\RR}, \  S_{t+1}=S_t$, preserving a smooth volume 
$\mu$. For any Liouville number $\a$ we construct a
sequence of area-preserving diffeomorphisms 
$H_n$  such that the sequence $H_n\circ S_\a\circ H_n^{-1}$
converges to a smooth weak mixing diffeomorphism of $M$. The method is a 
quantitative version of the approximation by conjugations
construction introduced in \cite{AK}. 

 For $m=2$ and  $M$ equal to the unit disc $\DD^2=\{x^2+y^2\leq 1\}$
or the closed annulus $\AAA=\TT\times [0,1]$  this result proves the following dichotomy: 
$\a \in \RR \setminus\QQ$ is Diophantine if and only if
there is no ergodic diffeomorphism of $M$ 
whose rotation number on the boundary equals $\alpha$ (on at least one of the boundaries in the case of 
$\AAA$). 
One part of the dichotomy follows from our constructions, the other is
an unpublished result of Michael
Herman asserting that if $\a$ is Diophantine, then any area preserving
diffeomorphism
with rotation number $\a$ on the boundary (on at least one of the boundaries in the case of $\AAA$) displays  
smooth
invariant curves arbitrarily close to the boundary which clearly
precludes ergodicity or even topological transitivity.  

\end{abstract}

\section{Introduction}
We present a construction method providing analytic weak mixing
diffeomorphisms on the torus $\TT^d=\RR^d/\ZZ^d$, $d\geq 2$, and smooth  
weak mixing diffeomorphisms on any smooth manifold with a nontrivial circle action preserving a smooth 
volume $\mu$. The diffeomorphisms obtained are homotopic to the Identity and
can be made arbitrarily close to it.

We will effectively work either on the two torus for the analytic constructions 
or on the closed annulus 
 $\AAA=\TT\times [0,1]$ for the smooth constructions.
In the case of the torus the construction is exactly the same in higher dimensions and we explain in \S 
\ref{reduction}  how the smooth construction can be transfered from the annulus to  general manifolds with a 
nontrivial circle action.

By smooth diffeomorphisms on a manifold with boundary we mean 
infinitely smooth in the interior and such that all the derivatives 
can  be continuously extended to the boundary. 

We recall that a dynamical system $(M,T, \mu)$ is said to be {ergodic} if and only if  there is no 
nonconstant {\it invariant } measurable complex function $h$ on $(M,\mu)$, i.e.  such that $h(Tx) = h(x)$. 
It is said to be weak mixing 
if  it enjoys the stronger property  of not having {\it eigenfunctions} at all, i.e. if 
 there is no nonconstant measurable complex function $h$ on $(M,\mu)$  such that $h(Tx) = \lambda h(x)$ for 
some constant $\lambda \in \CC$.

The construction, on any smooth manifold with a nontrivial circle action (in particular $\DD^2$), of volume 
preserving diffeomorphisms  enjoying different ergodic properties
(among  others, weak mixing)  was first undertaken in \cite{AK}. 
For $t \in \RR$ denote by $S_t$ the elements of the circle action on
$M$ with the normalization $S_{t+1}=S_t$. 

Let $\A(M)$ be the  closure  in the $C^\infty$ topology of the set of
diffeomorphisms of the form $h\circ S_t\circ h^{-1}$, with $t\in \RR$ and $h$
area preserving $C^\infty$-diffeomorphism of $M$.

 For a given $\a \in \RR$ we denote
by $\A_\a(M)$ the restricted space 
of {\sl conjugacies 
of the fixed rotation} $S_\a$, 
namely the closure of the set of  $C^\infty$-diffeomorphisms of
the form $h\circ S_\a\circ h^{-1}$. 

It is easy to see that the sets $\A_\a(M)$ are disjoint for different 
$\a$ and in \cite[section 2.3.1]{FK}, it was proved  for a particular 
manifold $M$ that $\cup_{\a \in \RR}  \A_\a(M) 
\varsubsetneq \A(M)$. We do not know if the inclusion remains strict 
on any manifold.

Anosov and Katok  proved in  \cite{AK} that  in $\A(M)$ the set of weak
mixing diffeomorphisms is generic (contains a $G_\delta$ dense set)
in the $C^\infty$ topology. 
Actually, it also
follows from the same paper that the same is true in $\A_\a(M)$ for
a $G_\delta$ dense set of $\a \in \RR $ although the construction,
properly speaking, is achieved in the space $\A(M)$. 
However, \cite{AK} does not give a full description of the set of $\a$
for which the result holds in $\A_\a(M)$. Indeed, the flexibility of
the constructions in \cite{AK} comes from the fact that $\a$ is
constructed inductively at the same time as the conjugations are
built, that is: at step $n$, $\a_n=p_n/q_n$ is given, and $h_n$ is
constructed that commutes with 
$S_{\a_n}$; then $\a_{n+1}$ is chosen so close to $\a_{n}$ that 
$f_n = H_n S_{\a_{n+1}}
H_n^{-1}$  (where $H_n = h_1 \circ \dots  \circ h_n$ and each $h_n$ commutes with $S_{\a_n}$) is 
sufficiently
close to $f_{n-1}$ to guarantee the convergence of the sequence  
${\lbrace f_n \rbrace}_{n \in \NN}$. Then
step $n+1$ gets started by the choice of $h_{n+1}$ etc... The final $\a$ is the limit of $\a_n$. 
By this
procedure, there is no need to put any restrictions on the growth of
the $C^r$ norms of $H_n$ since $\a_{n+1}$ can always be chosen close
enough to $\a_n$ to force convergence. The counterpart is that the
limit diffeomorphism obtained in this way will lie in $\A_\a(M)$ with
$\a$ having rational approximations at a speed that is not controlled.

Since we want to do the construction inside $\A_\a(M)$ for an
arbitrary Liouville number $\a$, 
we are only allowed to make use of the fact
that the decay of $|\a_{n+1} -
\a_n|$ is faster than any polynomial in $q_n$. So we have to 
construct
$h_n$ with a polynomial (in $q_n$) control on the  growth of its derivatives 
to make
sure that the above procedure converges.

Recall that an irrational number $\a$ is said to be Diophantine if it
is not too well approximated by rationals, namely if there exist
strictly positive constants $\gamma$ and $\tau$ such that for any
couple of integers $(p,q)$ we have:
$$ 
|q\a-p| \geq \frac{\gamma}{ q^{\tau}}.
$$

In this paper we
work in the restricted spaces  ${\mathcal A}_\a(M)$ and prove the
following for {\it any} Liouville, i.e. not Diophantine and not
rational,  frequency $\a$:
\begin{theo} \label{liouville}  Let $M$ be an $m$-dimensional ($m \geq 2$)
differentiable manifold with a nontrivial circle action
  ${\mathcal S} = {\lbrace S_t \rbrace }_{t \in\RR}, \  S_{t+1}=S_t$ preserving a smooth volume 
$\mu$. If $\a \in \RR $ is Liouville,
then the set of weak mixing diffeomorphisms  is generic 
in the $C^\infty$ topology in ${\mathcal A}_\a(M)$. 
\end{theo}

On $M=\DD^2$ or  $\AAA$, the weak mixing diffeomorphisms we will construct in 
${\mathcal A}_\a(M)$ will have $S_\a$ as their restriction to the boundary. This clarifies the relation 
between   
the ergodic properties of  the area preserving diffeomorphisms of $\DD^2$ and their 
rotation number on the boundary, complementing the
 striking result
of M. Herman stating that  if $f$ is a smooth diffeomorphism of the disc with a Diophantine rotation number 
on the boundary, then there exists a set of
positive measure of smooth invariant curves in the neighborhood of
the boundary, thus $f$ is not ergodic.  By KAM theory, this phenomenon was known to happen for
Diophantine $\a$ as soon as the map $f$ displays some twist features
near the boundary. Herman's tour de force  was to get rid of the
twist condition in the area preserving context. 
To be more precise, we introduce the following

\begin{defi}  Let $M$ denote either $\DD^2$ or  $\AAA$. 
Given $\a \in \RR$, we denote by ${\mathcal B}_\a(M)$
the set of area preserving $C^\infty$-diffeomorphisms of  $M$  whose
restriction to the boundary (to at least one of the boundary circles in the case of the annulus)  has a 
rotation number $\alpha$.
\end{defi}

\begin{theo}[Herman] \label{diophantien} Let $M$ denote 
either $\DD^2$ or  $\AAA$.
For a Diophantine $\a$, let $F\in {\mathcal B}_\a(M)$. Then the boundary
of $M$ (on which the rotation number is $\a$)  is accumulated by a set of positive measure of invariant 
curves of $F$.
\end{theo}
In the case of the disc and the annulus, as a corollary of Theorems 
\ref{liouville} and \ref{diophantien}, we have
 the following characterization of Diophantine numbers:
\begin{coro} \label{dichotomie}
Let $M$ denote either $\DD^2$ or  $\AAA$. A number 
$\a \in \RR \setminus \QQ$  is Diophantine if
and only if there is no ergodic diffeomorphism  $f \in {\mathcal B}_\a(M)$.  
\end{coro}

On $M= \TT^2$ and under a more restrictive condition on $\a$, the method
of  approximation by conjugations can be undertaken in the real
analytic topology and with very explicit conjugations. For an arbitrary fixed
$\si >0$, for any $n \in \NN$, we set:
\begin{equation}\label{eqdef_apprintro-1}
\begin{aligned} 
\phi_n (\theta,r) &= (\theta, r+q_n^2 \cos (2\pi q_n \theta) ), \\
g_n(\theta,r)&=(\theta+[n q_n^{\si}] r, r), \\
h_n&= g_n \ci \phi_n, \quad H_n =h_1 \circ \dots \circ  h_n, \\
f_n&=H_n \ci R_{\a_{n+1}} \ci H_n^{-1}.
\end{aligned}
\end{equation}
Here [$\,\cdot\,$] denotes the integer part of the number and $R_t$ denote the action $(\theta,r) \to 
(\theta+t,r)$.
The convergence of the diffeomorphisms $f_n$ is in the sense 
of a usual metric $d_\rho(\cdot,\cdot)$,  based on the supremum norm
of analytic functions over the complex strip of width $\rho$; see
Section \ref{analtop} for the definition.
We will prove the following 
\begin{theo} \label{analytiques}
Let $\a\in \RR$ be such that,  for some $\de>0$, equation
$$
|\a-p_n/q_n|<\exp (-q_n^{1+\de})
$$
has an 
infinite number of integer solutions $p_n$, $q_n$ (where 
$p_n$ and $q_n$ are relatively prime for each $n$). 
Take $\, 0<\si< \min\{\de/3, 1\}$. Then, for all $\rho>0$, 
there exists a sequence $\a_n=p_n/q_n$ (which is a subsequence of the 
solutions of the equation above) such that  
the corresponding diffeomorphisms $f_n$, constructed in 
$(\ref{eqdef_apprintro-1})$,
converge in the sense of the $d_\rho(\cdot,\cdot)$-metric,
and $f=\lim_{n\to\infty} f_n$ is weak mixing.
\end{theo}
Weak mixing diffeomorphisms, given by this theorem, are 
uniquely ergodic. This can be shown by the same method as in
\cite{S}. 

\begin{rema} The result in Theorem \ref{analytiques} is actually weaker
than what can be obtained by time change, e.g. the existence on
$\TT^2$  of real analytic weak mixing reparametrizations of
$R_{t(1,\a)}$ for  any irrational $\a$ such that $\mathop{\limsup}
\limits_{p \in \ZZ,q \in \NN^*} - \frac{ \ln|\a - {p / q}|}{ q}
\neq 0$ \cite{weakmixing,kr,kolmogorov,shklover}. Indeed, such reparametrizations  belong
{\sl a priori} to ${\mathcal A}_\a(\TT^2)$ (Cf. \cite{FK}). However, we
included the constructions on $\TT^2$ with explicit successive
conjugations as in (\ref{eqdef_apprintro-1}) because the proof of
weak mixing follows almost immediately from the  general criteria
we established to treat the general smooth case, and also
because these constructions might be generalized to other
manifolds where the techniques of reparametrizations are not
available. \end{rema}

{\it Acknowledgments.} It is a pleasure to thank H{\aa}kan Eliasson 
for useful comments all along the work. We are grateful to Rapha\"el
Krikorian for his  help in the main construction of Section \ref{section-smooth} 
and to Anatole Katok, Jean-Paul Thouvenot and Alistair Windsor 
for many useful conversations. We also thank the referee for many useful recommendations. 

\section{Preliminaries}

\subsection{General scheme of the constructions.}
 
Here we give a general scheme of the construction of the diffeomorphisms as a limit of conjugacies of a 
given Liouvillean action while \S \ref{cwm} outlines the particular choices that will yield the weak mixing 
property for the limit diffeomorphism. Henceforth, $M$ denotes either the torus $\TT^2$ or the annulus 
$\AAA$ and we consider polar coordinates $(\theta,r)$ on $M$ that denotes either the torus $\TT^2$ or the annulus 
$\AAA$. By $\la$ and $\mu$ we denote the usual Lebesgue measures on $\RR$ and
on $\RR^2$, respectively. The term ``measure-preserving" will refer to
the measure $\mu$.
 
For $\a \in \RR$, we consider the map $S_\a : M \to M$,
$(\theta,r) \mapsto (\theta+\alpha, r)$. The diffeomorphisms 
that we shall
construct, are obtained as limits  of measure preserving transformations 
\begin{eqnarray} \label{conjugations} 
f = \lim_{n \rightarrow \infty} f_n,\,\,\,\,\text{where}  
\,\,\,\, f_n=H_n\circ S_{\a_{n+1}}\circ H_n^{-1} .
\end{eqnarray} 
Here
$\a_n= p_n/q_n $ is a  convergent sequence of
rational numbers, such that $|\a-\a_n|\to 0$ monotonically;
$H_n$ is a sequence of measure preserving  diffeomorphisms of
$M$. In different constructions, the convergence  of $f_n$ will be
meant in the
$C^\infty$ or real analytic category; the topology in each case is
standard, and will be recalled in Sections \ref{analtop} and \ref{smoothtop}.

Each $H_n$ is obtained as a composition 
\begin{eqnarray} \label{conjugations2} 
H_n =  h_1 \circ \dots \circ h_n,
\end{eqnarray}
where every $h_n$ is a  measure preserving diffeomorphism of $M$ 
satisfying 
\begin{eqnarray} \label{conjugations3} 
h_n \circ S_{\a_n} = S_{\a_n} \circ h_n.
\end{eqnarray} 
At step $n$, $h_n$
must display enough stretching to insure an increasing distribution of
the orbits of $H_n\circ S_{\a_{n+1}}\circ H_n^{-1}$. 
However, this stretching must be
appropriately controlled  with respect to $|\a- \a_n|$ to guarantee
convergence of the construction.

\subsubsection{Decomposition of $h_n$}

In the subsequent constructions, each $h_n$ will be obtained as a composition 
\begin{eqnarray} 
\label{gn} h_n = g_n \circ \phi_n, 
\end{eqnarray}
where $\phi_n$ is constructed 
in such a way that 
$S_{1 / q_n} \circ \phi_n = \phi_n\circ S_{1 / q_n}$;
the diffeomorphism $g_n $ is a twist map  of the form 
\begin{eqnarray} \label{phin} 
g_n(\theta,r) = (\theta + [n q_n^{\sigma}] r , r), 
\end{eqnarray}
for some $0< \sigma < 1$ that will be fixed later.
The role of $g_n$ is to introduce shear in the "horizontal" direction
(the direction of the circle action), while $\phi_n$ is responsible
for the "vertical" motion, i.e. transversal to the circle action. 
The choice of the shear factor $nq_n^{\sigma}$ will be explained in \S \ref{overview}.

In the real analytic case, $\phi_n$ will be  given by an explicit
formula and convergence will follow from an assumption on the rational
approximations of $\a$. In the smooth case, $\phi_n$ will be
constructed in Section \ref{sectionphi} in such a way that its
derivatives satisfy   estimates of the type: 
$$
\|D_a\phi_n\|_{0} \leq c(n,a) q_n^{|a|},\quad
\|D_a\phi_n^{-1}\|_{0}\leq c(n,a) q_n^{|a|},
$$
where $c(n,a)$ is independent of $q_n$  (Cf. \S \ref{smoothtop} and \S \ref{notation} about
the notations we adopt). This polynomial growth 
of the norms of $\phi_n$ is crucial to insure the convergence of the
construction above and is the reason why it can be carried out for an
arbitrary Liouville number.

\subsection{Analytic topology.}\label{analtop}
Let us discuss the topology on the space of real-analytic diffeomorphisms of
$\TT^2$, homotopic to the identity. All of them have a
lift of type 
$F(\theta,r)=(\theta+f_1(\theta,r), r+f_2(\theta,r))$, 
where $f_i:\RR^2\to \RR$ are real-analytic 
and $\ZZ^2$-periodic. 

For any $\rho >0$, consider the set of real
analytic $\ZZ^2$-periodic functions on $\RR^2$, that can be
extended to holomorphic functions on $A^\rho=\{|\text{Im}\, \theta|, 
|\text{Im}\, r|<\rho\}$. For a function $f$ in this set, let
$\|f\|_\rho=\sup_{A^\rho}|f(\theta,r)|$.
We define $C^{\omega}_\rho(\,\TT^2)$ as a subset of the above set, 
defined by the condition: $\|f\|_\rho<\infty$.

Consider the space $\diff_\rho^\omega$ of those
diffeomorphisms, for whose lift it holds: $f_i\in
C^{\omega}_\rho (\,\TT^2)$, $i=1,2$.
For any two diffeomorphisms  $F$ and $G$ in this space 
we can define the distance
$$
d_\rho(F,G)=\max_{i=1,2}\{\,\inf_{p\in \ZZ} \|f_i-g_i+p\|_\rho\}.
$$

For a diffeomorphism $T$ with a lift  $T(\theta,r)=(T_1(\theta,r),\, T_2(\theta,r))$ 
denote 
$$
\|DT\|_\rho = \max \left\{ \left\|\frac{\partial T_1}{\partial \theta}\right\|_\rho, 
\left\|\frac{\partial T_1}{\partial r}\right\|_\rho,  
\left\|\frac{\partial T_2}{\partial \theta }\right\|_\rho,  
\left\|\frac{\partial T_2}{\partial r}\right\|_\rho
\right\}.
$$

\subsection{$C^\infty$-topology.}\label{smoothtop}

Here we discuss the (standard) topology 
on the space of smooth diffeomorphisms
of $M=\TT^2$, which we shall use later. 
The annulus is endowed with the topology in the similar way.

We are interested in convergence in the space 
of smooth diffeomorphisms of $M$, homotopic to the identity, 
and hence having lift of type $\tilde F(\theta,r)=(\theta+f_1(\theta,r),\, r+f_2(\theta,r))$,
where $f_i:\RR^2\to \RR$ are  $\ZZ^2$-periodic.
For a continuous function $f:(0,1)\times (0,1)\to \RR$, denote 
$$
\|f\|_{0}:=\sup_{z\in (0,1)\times (0,1)} |f(z)|.
$$
For conciseness we introduce the following notation for
partial derivatives of a function: for $ a=(a_1,a_2) \in \NN^2$ 
we denote $|a| := a_1 + a_2$ and
$$
D_a  := \frac{\partial^{a}}{\partial r^{a_1} \partial
\theta^{a_2} }.
$$

For $F$, $G$ in the space 
$\diff^k(\TT^2)$ of $k$-smooth diffeomorphisms of the torus, 
let  $\tilde F$
and $\tilde G$ be their lifts.
For mappings $F:\RR^2 \rightarrow \RR^2$ denote by 
$F_i$ the $i$-th coordinate
function.
Define the distances between two diffeomorphisms $F$ and  $G$ as
$$
\begin{aligned}
\tilde d_0(F,G)& = \max_{i=1,2} 
\{ \inf_{p\in \ZZ} 
\| (\tilde F -\tilde G)_i +p\|_{0}
\}, \\
\tilde d_k(F,G)& = \max \{ \tilde d_0(F,G), \ \| D_a (\tilde F_i -\tilde G_i)
\|_{0}
\mid 
i=1,2,\ \  1 \leq
|a| \leq k \}.
\end{aligned}
$$
We shall use the metric, measuring the distance both between
diffeomorphisms and their inverses:
$$
d_k(F,G)= \max \{\tilde d_k(F,G), \tilde d_k(F^{-1},G^{-1}) \}.
$$  


For $M=\DD^2$, the $\diff^k(M)$ topologies are defined in the natural 
way with the help of the supremum norm of continuous functions over
the disc.  

For the smooth topology on $M$, a sequence of $\diff^\infty(M)$ 
diffeomorphisms 
is said to be convergent in $\diff^\infty(M)$, if it converges 
in $\diff^k(M)$ for all $k$.
The space $\diff^\infty(M)$, endowed with the metric 
$$
d_\infty(F,G)=\sum\limits_{k=1}^\infty\frac{d_k(F,G)}{2^k(1+ d_k(F,G))},
$$
is a compact metric space, hence for any of its closed subspaces,
Baire theorem holds.

\subsection{Reduction to the case of the annulus} \label{reduction}

Let $(M,{\mathcal S},\mu)$ denote a system of an $m$-dimensional smooth 
manifold with a nontrivial circle action preserving a smooth volume $\mu$. 

We denote by $F$ the set of  fixed points of the action ${\mathcal S}$.
For $q \geq 1$ we denote by $F_q$ the set of fixed points of the map $S_{1/q}$.
And by $\partial M$ we denote the boundary of $M$. Finally we let $B:= \partial M \cup F \cup_{q\geq1} F_q$. 

Let $\lambda$ be the product of Lebeasgue measures on $\LL^1 \times \DD^{m-1}$.  
Denote by ${\mathcal R}$ the standard "horizontal" action of $\LL^1$ on  $\LL^1 \times \DD^{m-1}$.  We quote 
the following proposition  of \cite{FK} that is similar to corresponding 
statements in \cite{AK,windsor}

\begin{prop}\cite[proposition 5.2]{FK} \label{reductiondirect} Let $M$ be an $m$-dimensional  differentiable 
manifold with a non-trivial circle action
  ${\mathcal S} = {\lbrace S_t \rbrace }_{t \in\RR}, \  S_{t+1}=S_t$ preserving a smooth volume 
$\mu$.  Let $B:=\partial M\cup F\cup(\underset{q}\bigcup F_q)$.  There exists a continuous surjective map 
$G: \LL^1 \times {\DD}^{m-1} \to M$  
with the following properties:
\begin{enumerate}
\item The restriction of $G$ to the interior $ \LL^1 \times  {\DD}^{m-1}$  
is a $C^{\infty}$ diffeomorphic embedding;
\item $\mu(G(\partial(\LL^1 \times {\DD}^{m-1})))=0$;
\item $G(\partial(\LL^1 \times {\DD}^{m-1})) \supset B$;
\item$G_*(\lambda)=\mu$;
\item $\mS\circ G=G\circ\mathcal R$.
\end{enumerate}
\end{prop}
We show now how this proposition allows to carry a construction as in the preceding 
section from $(\LL^1 \times \DD^{m-1}, {\mathcal R}, \lambda)$ to the general 
case $(M,{\mathcal S},\mu)$.

Suppose $f: \LL^1 \times {\DD}^{m-1} \to \LL^1 \times  {\DD}^{m-1}$ is a weak mixing 
diffeomorphism given, as above, by $f= \lim f_n$, $f_n=H_n  \circ R_{\a}  \circ H_n^{-1}$ 
where, moreover, the maps $H_n$ are equal to identity in a neighborhood of the boundary, 
the size of which can be chosen to decay arbitrarily slowly. Then if we define the 
diffeomorphisms $\tilde{H}_n: M \to M$ 
\begin{eqnarray*}
\tilde{H}_n(x) &=& G \circ H_n \circ G^{-1} (x) 
\ {\rm for \ } x \in G ( \LL^1 \times  {\DD}^{m-1}), \ {\rm and} \\
\tilde{H}_n(x) &=& x \ {\rm for \   } x \in  G ( \LL^1 \times \partial ({\DD}^{m-1})),
\end{eqnarray*}
we will have that $\tilde{H}_n  \circ S_\a  \circ \tilde{H}_n^{-1} $ is convergent in the $C^\infty$ 
topology to the weak mixing diffeomorphism
$\tilde{f} : M \rightarrow M$ defined by
$$g(x)=G(f(G^{-1}(x)) \ {\rm for \ } x \in G ( \LL^1 \times  {\DD}^{m-1}), \ {\rm and}$$
$$g(x) = S_\a(x) \ {\rm for \   } x \in  G ( \LL^1 \times  \partial ({\DD}^{m-1})).$$

In the sequel, to alleviate the notations, we will assume that $m=2$ and will do the constructions on the 
annulus $\AAA=  \LL^1 \times  [0,1]$ or on the two torus $\TT^2$.

\section{Criterion for weak mixing}\label{cwm}

The goal of this section is to give a simple geometrical criterion  involving only the
diffeomorphisms  $\phi_n \circ R_{\alpha_{n+1}}\circ \phi_n^{-1}$ and insuring the weak mixing property for 
the diffeomorphism $f$
given by (\ref{conjugations})--(\ref{phin}) in case of
convergence. The criterion will be stated in Proposition
\ref{criterion} of \S \ref{cccc}.

The following characterization of weak mixing will be used 
(see, for example, \cite{shklover}): $f$ is weak
mixing if there exists a sequence $m_n \in \NN$ such that for any
Borel sets $A$ and $B$ we have: 
\begin{eqnarray} \label{wwww}
|\mu(B\cap f^{-m_n}(A))-\mu(B)\mu(A)|\to 0.
\end{eqnarray}

\subsection{} \label{overview} We will now give an overview of the criterion assuming that $M$ is the
annulus $\TT \times [0,1]$ and denoting by horizontal intervals the
sets $I= [\theta_1,\theta_2] \times \lbrace r \rbrace$. We say that a
sequence $\nu_n$, consisting for each $n$ of a collection of disjoint
sets on $M$ (for example horizontal intervals), {\it converges to the 
decomposition into points} if any measurable set $B$  can be
approximated as $n \rightarrow \infty$ by a union of atoms in $\nu_n$
(Cf. \S \ref{ffff}). We denote this by $\nu_n \mathop{\longrightarrow}
\limits_{n \rightarrow \infty} \epsilon$.

The first reduction  is given by a  Fubini Lemma
\ref{Fubinilemma}. Here we 
decompose $B$  at each step $n$ into a union of small codimension
one sets for  which a precise  version of (\ref{wwww}) is assumed to hold, see (\ref{fub1}). 
For each $n$ these sets are images by a smooth map $F_n$ of a collection $\eta_n$ of
horizontal intervals such that $F_n(\eta_n) \rightarrow \epsilon$.
Lemma \ref{Fubinilemma} shows that (\ref{fub1}) guarantees weak mixing.

The second step is Lemma \ref{re2} 
asserting
that under an additional condition of proximity  (\ref{zvezdochka})
between $f_n^{m_n}$ and $f^{m_n}$, it is enough to  check (\ref{fub1})
for $f_n$.

Now, we take $F_n$ in the Fubini Lemma equal to $H_{n-1} \circ g_n$. 
Since $H_{n-1}$ in the construction
only depends on $q_{n-1}$, $q_n$ can be chosen so that $\|DH_{n-1}\|_{0}<\ln q_n$. With our choice 
of $g_n$ ($\sigma <1$ in (\ref{phin})) this implies that 
$ H_{n-1} \circ g_n (\eta_n) \rightarrow \epsilon$ 
if $\eta_n \rightarrow \epsilon$ is a partial partition with horizontal inetrvals of length less than 
$1/q_n$ (Cf. Lemma \ref{nu->e}). With the above observations, we are reduced to finding a collection 
$\eta_n$ and a sequence $m_n$ with the property that $H_{n-1} \circ g_n \circ \phi_n \circ 
R_{\a_{n+1}}^{m_n} \circ \phi_n^{-1} (I)$ is almost uniformly distributed in $M$ for $I \in \eta_n$.  
  
The geometrical ingredient of the criterion appears in \S \ref{hhhh} and
merely states that if a set (in particular $\phi_n \circ R_{\a_{n+1}}^{m_n} \circ \phi_n^{-1} (I)$) is 
almost a vertical line going from one
boundary of the annulus to the other, then the image of this set by
$g_n$ defined in (\ref{phin}) is  almost uniformly distributed in
$M$. "Almost vertical"  is made precise and quantified in Definition
\ref{aaaa}. Actually, the choice of $g_n$ ($\sigma > 0$ in \ref{phin}) gives in addition that 
$H_{n-1} \circ g_n$ of a an almost vertical segment will be almost  uniformly distributed in
$M$, since we impose that $\|DH_{n-1}\|_{0}<\ln q_n$.

In conclusion, the criterion for weak mixing (Proposition \ref{criterion}) roughly states as follows: 
Let $f$
be given by  (\ref{conjugations})--(\ref{phin}). If for some sequence
$m_n$   satisfying the proximity condition (\ref{zvezdochka})  between
$f_n^{m_n}$ and $f^{m_n}$, there exists a sequence $\eta_n \rightarrow
\epsilon$ consisting of horizontal intervals of length less than
$1/q_n$ such that the image of $I \in \eta_n$ by $\phi_n \circ R_{\a_{n+1}}
\circ \phi_n^{-1}$ is  increasingly almost vertical as $n \rightarrow
\infty$ then the limit diffeomorphism $f$ is weak mixing.


\subsection{A Fubini Lemma.} \label{ffff}

\begin{defi}
A collection of disjoint sets on $M$ 
will be called {\it partial decomposition} of $M$.
We say that a sequence of  
partial decompositions $\nu_n$
converges to the decomposition into points (notation: $\nu_n\to \epsilon$)
if, given a measurable set $A$, 
for any $n$ there exists a measurable set $A_n$,
which is a union of elements of $\nu_n$, such that $\lim_{n\to \infty} 
\mu( A \triangle
A_n) =0$ (here $\triangle$ denotes the symmetric difference).
\end{defi}
In this section we work with 
$M= \TT^2$ or $M= \AAA$. For these manifolds we formulate the
following definition.
\begin{defi}
Let $\hat \eta$ be a  partial decomposition of
$\TT$ into intervals, and consider on 
$M$ the decomposition $\eta$ consisting of intervals in $\hat{\eta}$ times some $r \in [0,1]$.
Decompositions of the above type will be called {\it standard
partial decompositions}. We shall say that $\nu$ is the image under 
a diffeomorphism $F:M\to M$ of a 
standard decomposition $\eta$ (notation: $\nu=F(\eta)$), if 
$$
\nu=\{\Gamma=F(I)\mid I \in \eta \}. 
$$
\end{defi}

Here we formulate a standard criterion for weak mixing. The proof is
based on the application of Fubini Lemma.

\begin{lemm}[Fubini Lemma]\label{Fubinilemma} 
Let $f$ be a measure $\mu$ preserving diffeomorphism of $M$. 
Suppose that there exists an increasing sequence $m_n$ of 
natural numbers, and a sequence of partial   
decompositions $\nu_n\to \epsilon$ of $M$, where, for each $n$,  
$\nu_n$ is the 
image under a measure-preserving diffeomorphism $F_n: M\to M$ of a  
standard partial decomposition, with the following property: for any
fixed square $A \subset M$ and any $\varepsilon > 0$, for any $n$
large enough we have: for any atom 
$\Gamma_n \in \nu_n$ 
\begin{eqnarray}\label{fub1} 
\left|\la_n (\Gamma_n \bigcap f^{-m_n} (A) ) - 
\la_n(\Gamma_n) \mu(A)  \right| 
\leq \varepsilon \la_n(\Gamma_n)\mu(A),
\end{eqnarray}
where $\la_n =F_n^*(\la )$.

Then the diffeomorphism $f$ is weak mixing.
\end{lemm}

\begin{proof}
To prove that $f$ is weak mixing, 
it is enough to show that for any  square $A$ and a Borel set $B$
$$
|\mu(B\cap f^{-m_n}(A))-\mu(B)\mu(A)|\to 0
$$
when $n\to \infty$. In the case of the annulus, is even enough to 
show this  for any square $A$ that is strictly contained in the
interior of $\AAA$. 
By assumption, for any $n$ we have: 
$\la_n(\Gamma _n)=\la_n(F_n (I_n))=\la(I_n)$. Then
$$
\la_n(\Gamma_n\cap f^{-m_n} (A))=\la_n(F_n( I_n \cap F_n^{-1}\circ f^{-m_n} (A))
)=\la(I_n \cap F_n^{-1}\circ f^{-m_n} (A)).
$$
By (\ref{fub1}), this implies: 
$$
|\la(I_n \cap  F_n^{-1}\circ f^{-m_n} (A))- \la(I_n)\mu(A) |\leq \e \la(I_n)\mu(A).
$$

Take any Borel set $B\subset \TT^2$. Since $\nu_n\to \epsilon$, for
any $\e$, for fixed $A$ and $B$, there exists $n$ and a measurable set 
$\hat B=\cup_{i\in \si} \Gamma_n^i $ ($\Gamma_n^i$ are elements of $\nu_n$,
and $\si$ is an appropriate index set) 
such that 
$$
|\mu(\hat B \triangle B)|<\e\mu(B)\mu(A).
$$
Consider $\tilde B=F_n^{-1}(\hat B)$ (it is also measurable since $F_n$
is continuous). 
Then 
$$
\tilde B=\bigcup_{i\in \si}  F_n^{-1}( \Gamma_n^i )=\bigcup_{i\in \si}
I_n^i:=\bigcup_{0\leq y \leq 1} \bigcup_{i\in\si(y)} I^i_n(y)  \times \{ y \}.
$$
We estimate:
$$
\begin{aligned}
&|\mu(B\cap f^{-m_n}(A))-\mu(B)\mu(A)| \\
&= |\mu(F_n^{-1}(B)\cap F_n ^{-1} \circ
f^{-m_n}(A))-\mu(B)\mu(A)|  \\
&\leq |\mu(\tilde B \cap F_n ^{-1}\circ  
f^{-m_n}(A))-\mu(\tilde B)\mu(A)|+2\e\mu(B)\mu(A)  \\
&=\int_{0}^{1} \sum_{i\in\si(y)} 
| \la(I^i_n(y)  \times \{ y \} \cap F_n^{-1}\circ  f_n^{-m_n} (A) ) -
\la(I^i_n) \mu (A) |  \text{d}y \\
&+2\e\mu(B)\mu(A) \leq 
 3 \e \mu(B) \mu(A).  
\end{aligned}
$$

\end{proof}
 
\subsection{Reduction from $f$ to $f_n$.} \label{rrrr}

\begin{lemm}[Reduction to $f_n$] \label{re2} 
If $f$ is the limit diffeomorphism from (\ref{conjugations}), and 
the sequence $m_n$ in the latter lemma satisfies 
\begin{equation}\label{zvezdochka} 
d_0(f^{m_n},f_n^{m_n})<\frac{1}{2^n},
\end{equation}
then
we can replace the diffeomorphism $f$ in the
criterion (\ref{fub1}) by $f_n$: 

\begin{eqnarray}\label{fub2} 
\left|\la_n (\Gamma_n \bigcap f_n^{-m_n} (A) ) - \la_n(\Gamma_n) \mu(A)  \right| 
\leq \varepsilon \la_n(\Gamma_n)\mu(A),
\end{eqnarray}
and the result of Lemma \ref{Fubinilemma} still holds.
\end{lemm}

\begin{proof}
Let us show that the assumptions of this lemma  imply (\ref{fub1}).
Fix an arbitrary square $A \subset M$ and $\varepsilon > 0$. 

Consider two squares, $A_1$ and $A_2$,
such that
$$
A_1\subset  A \subset A_2, 
\quad \mu(A\triangle A_i) \leq \frac{\e}{3} \mu(A).
$$ 
Moreover, if $n$ is sufficiently large, we can guarantee that
$$
\text{dist\,} (\partial A, \partial A_i)> \frac{1}{2^{n}},
$$
(where $\text{dist\,} (A,B) =\inf_{x\in A,\, y\in B}|x-y|$, and
$\partial A$ denotes the boundary of $A$), and
$$
\left| \la_n (\Gamma_n \bigcap f_n^{-m_n} (A_i) ) - \la_n (\Gamma_n) \mu(A_i)  \right| 
\leq \frac{\e}{3} \la_n(\Gamma_n)\mu (A_i).
$$

By (\ref{zvezdochka}), for any $x$ the following holds:
$f_n^{m_n}(x)\in A_1$ implies $f^{m_n}(x)\in A$, and 
$f^{m_n}(x)\in A$ implies $f_n^{m_n}(x)\in A_2$.
Therefore, 
$$
\la_n (\Gamma_n\cap f_n^{-m_n}(A_1))\leq \la_n (\Gamma_n\cap f^{-m_n}(A)) 
\leq \la_n (\Gamma_n\cap f_n^{-m_n}(A_2)),
$$
which gives the estimate:
$$
\left(1-\frac{\e}{3}\right)\la_n(\Gamma_n)\mu(A_1) 
\leq \la_n (\Gamma_n\cap f^{-m_n}(A)) 
\leq \left(1+\frac{\e}{3}\right) \la_n(\Gamma_n) \mu(A_2),
$$
implying (\ref{fub1}).
\end{proof}


\subsection{Reduction from $f_n$ to 
$h_n \circ R_{\a_{n+1}}\circ  h_n^{-1}$.} \label{nnnn}

The following is a technical lemma that will allow us to focus in the
sequel only on the action of $h_n\circ  R_{\a_{n+1}}^{m_n}\circ 
h_n^{-1}$ (more
specifically on $g_n \circ \phi_n \circ R_{\a_{n+1}}^{m_n} \circ
\phi_n^{-1}$)  in order to get (\ref{fub2}):

\begin{lemm}\label{nu->e}
Let $\eta_n$ be a sequence of standard partial decompositions of $M$ into
horizontal intervals of length less or equal to $1 / q_n$, let $g_n$ be
defined by $(\ref{phin})$ with some $0< \si<  1$, and  let $H_n$
be a sequence of area-preserving diffeomorphisms of $M$ such that
for all $n$ 
\begin{equation}\label{estimate}
\|DH_{n-1}\|_{0}<\ln q_n.  
\end{equation}
Consider partitions $\nu_n=\{\Gamma_n =H_{n-1} g_n (I_n)\mid I_n\in \eta_n
\}$.

Then $\eta_n\to \epsilon$ implies  $\nu_n\to \epsilon$.
\end{lemm}
\begin{proof}
Let $\sigma < \sigma' <1$, and consider a partition of the annulus into squares $S_{n,i}$ of 
side length between $q_n^{-\si'}$ and $2 q_n^{-\si'}$. Since $\eta_n \rightarrow \epsilon$, we have for 
$\varepsilon >0$ arbitrarilly small, if $n$  is large enough,
$\mu(\cup_{I \in \eta_n}I) \geq 1 - \varepsilon$, 
so that for a collection of atoms $S$ with total measure greater than $1 - \sqrt{\varepsilon}$ we have  
$\mu(\cup_{I\in \eta_n}I\cap S) \geq (1 -\sqrt{\varepsilon}) \mu(S)$. Since $\si'< 1$ and  any $I \in 
\eta_n$ 
has length at most 
$1 / q_n$, we have for the same atoms $S$ as above $\mu(\cup_{I\in \eta_n, I\subset S} ) \geq 
(1-2\sqrt{\varepsilon}) \mu(S)$ if $n$ is sufficiently large. 

Consider now the sets  $C_{n,i}= H_{n-1}g_n(S_{n,i})$. 
In the same way as the squares $S_{n,i}$, a large proportion of 
these sets can be  well approximated by  unions of elements of
$\nu_n$. But by (\ref{estimate}),  we have:
$$
\text{diam\,}(C_{n,i}) \leq \|DH_{n-1}\|_{0} \, \|Dg_{n}\|_{0} \, 
\text{diam\,}(S_{n,i}),
$$ 
which goes to $0$ as $n \rightarrow \infty$. Therefore, any Borel set $B$  
can be approximated by a union of such sets
$C_{n,i}$ with any ahead given accuracy, if $n$ is sufficiently large, 
hence $B$ gets well approximated by  unions 
of elements of $\nu_n$. 
\end{proof}

             
\subsection{Horizontal stretch under $g_n$.} \label{hhhh}

We shall call by {\it horizontal interval} any line segment of the
form $I\times \{r\} $, where $I$ is an interval on the $\theta$-axis.
{\it Vertical intervals} have the form $\{ \theta\} \times J$ where 
$J$ is an interval on the $r$-axis.
Let $\pi_r$ and $\pi_\theta$ denote the projection operators onto
$r$ and $\theta$ coordinate axes, respectively.

The following definition formalizes the notion of ``almost uniform
distribution" of a horizontal interval in the vertical direction.
 
\begin{defi}[$(\ga,\de,\e)$-distribution] \label{aaaa}
We say that a diffeomorphism $\Phi:M\to M$ $\ (\ga,\de,\e)$-distributes a
horizontal interval  
$I$ (or $\Phi(I)$ is $(\ga,\de,\e)$-distributed),
if 

$\bullet$ $\pi_r(\Phi(I) )$ is an interval $J$ with $1-\de\leq \la(J)\leq 1$;

$\bullet$  $\Phi(I)$ is contained in a ``vertical strip" 
of type $[c,c+\ga]\times J$ for some $c$;

$\bullet$ for any  interval $\ti J \subset J$ 
we have:
\begin{equation}\label{distribution}
\left|\frac{\lambda (I\cap \Phi^{-1} (\TT\times \ti J))}{\lambda(I) } -
\frac{\la(\ti J)}{\la (J)}\right|
\leq \e \frac{\la(\ti J)}{\la (J)}.
\end{equation}
\end{defi}
We shall more often write the latter relation in the form
$$
|\lambda (I\cap \Phi^{-1} (\TT\times \ti J)) \la(J)-\lambda(I)\la(\ti
J)|
\leq \e \lambda(I)\la(\ti J).
$$

\begin{lemm}\label{prelim2}
Let $g_n$ be a diffeomorphism of the form (\ref{phin}) with some fixed
$0<\si < 1$.
Suppose that a diffeomorphism $\Phi:M\to M$ $\ (\ga,\de,\e)$-distributes a
horizontal interval  
$I$ with $\ga=1/(nq_n^{\si})$, $\de=1/n$,
$\e=1/n$. 
Denote $\pi_r (\Phi(I))$ by $J$.

Then 
for any square $S$ of side length  $q_n^{-\si}$, lying in 
$\TT\times J$ it holds:
\begin{eqnarray}\label{eqprelim}
|\la (I \cap \Phi^{-1} \ci g_n^{-1} (S)) \la(J) -\la (I)\mu(S)|\leq
8/n \la(I)\mu(S).
\end{eqnarray} 
\end{lemm}

Lemma \ref{prelim2} asserts that, if a diffeomorphism $\Phi$
``almost uniformly" distributes $I$ in the vertical direction, then the
composition of  $\Phi$ and the affine map $g_n$ ``almost uniformly"
distributes $I$ on the whole of $M$.

To prove Lemma \ref{prelim2}, we shall need the following preliminary
statement: it says that $g_n$ ``almost uniformly" distributes on  $M$ any
sufficiently thin vertical strip.

\begin{lemm}\label{prelim1} Suppose that  $g:M\to M$ has a lift
$$
g(\theta, r)=(\theta + br, r)  \quad \text{for some }\ b\in \ZZ,\ 
|b|\geq 2.
$$    
For an interval $K$ on the $r$-axes, $\la(K)\leq 1$, denote by $K_{c,\ga}$ a strip
$$
K_{c,\ga}:= [c,c+\ga] \times K.
$$
Let $L=[l_1,l_2]$ be an interval on the $\theta$-axes. 
If $b\la(K)> 2$, then for
$$
Q:=\pi_r( K_{c,\ga} \cap g^{-1}(L\times K)),
$$
it holds: 
$$
|\la(Q)-\la(K)\la(L)|\leq  \ga \la(K)+ \frac{2\la(L)}{b}
  +\frac{ 2 \ga}{b}  . 
$$ 
\end{lemm}

\begin{proof}
By definition, $Q=\{ r \in K \mid \ \exists
  \theta\in [c,c+\gamma] : \theta+br \in [l_1,l_2] \} $.
Then
$$
Q =\{ r\in K \mid br \in [l_1-\gamma ,l_2]-c \}.
$$
To estimate $\la(Q)$, note that the interval $bK$ (seen as an interval
on the real line) intersects not more
than $b\la(K)+2$ intervals of type $[i,i+1]$, $i\in \ZZ$, on the line, and
not less than $b\la(K)-2$ such intervals. Hence,
$$
\la(Q)\leq
(b\la(K)+2)\frac{(l_2-l_1)+\gamma}{b}=\la(K)\la(L)+
\ga \la(K)+\frac{2 \la(L)}{b}+\frac{ 2 \ga}{b}.
$$
The lower bound is obtained in the same way.
\end{proof}

\begin{proof}[Proof of Lemma \ref{prelim2}] Let $S$ be a
square in $\TT \times J$ of size  $q_n^{-\si} \times q_n^{-\si}$. 
Denote $\pi_\theta(S)$ by $S_\theta$, 
$\pi_r (S)$ by $S_r$. In these notations, $\la(S_r)=\la(S_\theta)=
q_n^{-\si}$, and
$\la(S_\theta)\la(S_r)=\mu(S)=q_n^{-2\si}$.
  
Let us study what part of $\Phi(I)$ is sent by $g_n$ into $S$.
Since $\Phi(I)$ is contained in a strip $[c,c+\ga]\times
J$ for some $c$, by assumption, and $g_n$ preserves horizontals,
this part lies in $K_{c,\gamma}:=[c,c+\ga]\times S_r$. 
Denoting $S_\theta$ by $[s_1, s_2]$, define a ``smaller" rectangle
$S_1\subset S$: $S_1=[s_1+\ga,
s_2-\ga] \times S_r$ (in our assumptions, $2\ga$ is much less than
$\la(S_\theta)$, so this rectangle is non-empty). 
Consider two sets:
$$
Q:=\pi_r(K_{c,\gamma} \cap g_n^{-1}(S) ), 
\quad Q_1:=\pi_r(K_{c,\gamma} \cap g_n^{-1}(S_1) ).
$$
Then we have:
\begin{equation}\label{Oooo}
\Phi(I) \cap (\TT\times Q_1) \subset \Phi(I) \cap  g_n^{-1} (S)\subset \Phi(I) \cap (\TT\times Q).
\end{equation}
The second inclusion is evident, the first one comes from the fact that
$g_n$ preserves lengths of
horizontal intervals. 

Lemma \ref{prelim1} permits us to estimate $\la (Q)$ and $\la (Q_1)$.
Indeed, to estimate the former one, apply Lemma \ref{prelim1} with $b=[nq_n^{\si}]$,
$\ga=(nq_n^{\si})^{-1}$, $K=S_r$, and $L=S_\theta$. We get: 
$$
|\la(Q)-\mu(S)|\leq  \frac{\la(S_r)}{nq_n^{\si}}+ 
\frac{2\la(S_\theta)}{[nq_n^{\si}]}+ \frac{2}{nq_n^{\si}[nq_n^{\si}]} \leq \frac{4}{n} \mu(S) . 
$$
In the same way, applying  Lemma \ref{prelim1} with 
the same $b$, $\gamma$, $K$ as above and $L=\pi_\theta S_1= [s_1+\gamma,
s_2-\gamma] $,
we get the same estimate (for large $n$):
$$
|\la(Q_1)-\mu(S_1)|\leq\frac{4}{n} \mu(S).
$$
In particular, this implies $\la(Q)\leq 2\mu(S)$, and $\la(Q_1)\leq
2\mu(S)$.

Both $Q$ and $Q_1$ are finite unions of disjoint intervals.
Then, using (\ref{distribution}) with $\e=\frac{1}{n}$ (which was the
assumption of the present lemma), we have:
$$
|\la (I \cap \Phi^{-1}  (\TT\times Q) )\la(J) -\la (I)\la(Q)|
\leq \frac{1}{n}\la(I)\la(Q) \leq \frac{2}{n}\la(I)\mu(S),
$$
and the same estimate holds for $Q_1$ instead of $Q$.
The last preliminary estimates are:
$$
\begin{aligned}
&|\la (I \cap \Phi^{-1}  (\TT\times Q) )\la(J) -\la (I)\mu(S)| \leq \\
&|\la (I \cap \Phi^{-1}  (\TT\times Q) )\la(J) -\la (I)\la(Q)|
+\la(I)|\la(Q)-\mu(S)| \\
&\leq \frac{2}{n}\la(I)\mu(S)+\frac{4}{n}\la(I)\mu(S)=\frac{6}{n}\la(I)\mu(S); 
\end{aligned}
$$
and, in the same way, (noting that $\mu(S)-\mu(S_1)=
\frac{2}{n}\mu(S)$), 
one estimates 
$$
|\la (I \cap \Phi^{-1}  (\TT\times Q_1)
)\la(J) -\la (I)\mu(S)|\leq \frac{8}{n}\la(I)\mu(S).
$$ 

Now relation (\ref{Oooo}), together with  the preliminary estimates 
above, gives  the desired conclusion:
$$
\begin{aligned}
&|\la (I \cap \Phi^{-1} \ci g_n^{-1} (S)) \la(J)-\la (I)\mu(S)| \\
&\leq \max\{ |\la (I \cap \Phi^{-1}  (\TT\times Q) )\la(J) -\la (I)\mu(S)|,\\
&|\la (I \cap \Phi^{-1}  (\TT\times Q_1) )\la(J) -\la (I)\mu(S)|\}  
 \leq \frac{8}{n}\la(I)\mu(S).
\end{aligned}
$$
\end{proof}



\subsection{Criterion for weak mixing.} \label{cccc}

We can now state the following

\begin{prop}[Criterion for weak mixing] \label{criterion} 
Assume that $f_n = H_n \circ R_{\a_{n+1}}\circ  H_n^{-1}$ 
is a sequence of diffeomorphisms
constructed following 
(\ref{conjugations2}), (\ref{conjugations3}), (\ref{gn}) and
(\ref{phin}) with some $0<\si<1/2 $, and that for all $n$ (\ref{estimate}) holds.

Suppose that the limit $\lim\limits_{n\to \infty} f_n=f$ exists.
If there exist a sequence $m_n$ satisfying  
(\ref{zvezdochka}) and a
sequence of standard partial decompositions $\eta_n$ of $M$ into
horizontal intervals of length less than $1 / q_n$ such that 
\begin{enumerate} 
\item $\eta_n \rightarrow \epsilon$,

\item  for any interval $I_n \in \eta_n$, the diffeomorphism
$$
\Phi_n:=\phi_n \circ R_{\a_{n+1}}^{m_n}\circ  \phi_n^{-1} 
$$  
$(\frac{1}{nq_n^{\si}}, \frac{1}{n},\frac{1}{n})$-distributes the interval $I_n$, 
\end{enumerate}
then the limit diffeomorphism  $f$ is weak mixing. 
\end{prop}

\begin{proof}  We
use Lemma \ref{re2} to prove weak mixing.
Consider partitions 
$\nu_n=\{\Gamma_n =H_{n-1}\circ  g_n (I_n)\mid I_n\in \eta_n
\}$, and let $\la_n=(H_{n-1}\circ  g_n)^* \la $. 
By Lemma \ref{nu->e},    $\nu_n$  
converges to the decomposition into points.

Let an arbitrary square $A$ and $\e>0$ be fixed.
In order to be able to apply Lemma \ref{re2}, it is left to  
check condition (\ref{fub2}) for any $\Gamma_n \in \nu_n$, 
with $f_n^{m_n} = 
H_n \circ S^{m_n}_{\a_{n+1}}\circ  H_n^{-1}=
H_{n-1}\circ g_n \circ \Phi_n \circ g_n^{-1}\circ  H_{n-1}^{-1}$.
By assumption (2) of the present lemma, for all $I_n\in \eta_n$, 
$\pi_r(\Phi_n(I_n))\supset
[-1/n,1-1/n]$. Let $S_n$
be a square of side length
$q_n^{-\si}$, $S_n\subset \TT\times [-1/n,1-1/n]$. Consider 
$$ 
C_n:=H_{n-1}(S_n).
$$
Assumption  (2) permits to apply Lemma
\ref{prelim2}. Then we have (estimating $ \frac{1}{\la(J)}\leq 2$):
$$
\begin{aligned}
&|\la_n(\Gamma_n \cap f_n^{-m_n}(C_n))-\la_n(\Gamma_n)\mu(C_n)| \\
& = 
|\la(I_n \cap \Phi_n^{-1}\circ g_n^{-1}(S_n))-\la(I_n)\mu(S_n)| \\
&\leq \frac{1}{\la(J)}|\la(I_n \cap \Phi_n^{-1}\circ g_n^{-1}(S_n)) \la(J)
-\la(I_n)\mu(S_n)|
+ \frac{(1-\la(J))}{\la(J)}\la(I_n)\mu(S_n) \\
&\leq 2\frac{8}{n} \la(I_n)\mu(S_n)+\frac{2}{n} 
\la(I_n)\mu(S_n)=\frac{18}{n}\la_n(\Gamma_n)\mu(C_n).
\end{aligned}
$$

By (\ref{estimate}), we have for $n$ sufficiently large
$
\text{diam\,}(C_n) \leq  \|D(H_{n-1})\|_{0} \, \text{diam\,}(S_n)\leq \frac{1}{2^n}
$.
Hence, for $n$ large enough, one can
approximate $A$ by such sets $C_n$ lying in $\TT\times
[1/n,1+1/n]$. More precisely, for $n$ large enough,
there exist two sets, which are unions of sets $C_n$:
$A_1=\cup_{\sigma_1} C_n$, $A_2=\cup_{\sigma_2} C_n$ such that
$$
A_i \subset \TT\times [1/n,1-1/n], \quad 
A_1\subset A \cap \TT \times [1/n,1-1/n]  \subset A_2,
$$
$$
|\mu(A)-\mu(A_i)|\leq \frac{\e}{3}\mu(A).
$$
Take $n$ so that $\frac{18}{n}<\frac{\e}{3}$.
Then we can estimate:
$$
\begin{aligned}
&\la_n(\Gamma_n\cap f_n^{-m_n} (A))-\la_n(\Gamma_n)\mu(A)\leq
\la_n(\Gamma_n\cap f_n^{-m_n}
(A_2))-\la_n(\Gamma_n)\mu(A_2)+ \\
&\frac{\e}{3}\la_n(\Gamma_n) \mu(A) \leq 
\frac{\e}{3} \la_n(\Gamma_n)\mu(A_2)+
\frac{\e}{3}\la_n(\Gamma_n) \mu(A) \leq
\e\la_n(\Gamma_n) \mu(A).
\end{aligned}
$$
The lower estimate for this difference is obtained in the same way (using $A_1$).
We have shown that, if $n$ is sufficiently large,  
for an arbitrary $\Gamma_n\in\nu_n$  (\ref{fub2}) holds. 
Then, by Lemma \ref{re2}, $f$ is weak mixing.
\end{proof}





\section{Analytic case on the torus ${\TT}^2$.} \label{section-analytic}

  This section is devoted to the analytic construction on the torus $\TT^2$.
We recall the notations of the Theorem \ref{analytiques} that we want to prove. For an arbitrary fixed
$\si >0$, for any $n \in \NN$:
\begin{equation}\label{eqdef_apprintro}
\begin{aligned} 
\phi_n (\theta,r) &= (\theta, r+q_n^2 \cos (2\pi q_n \theta) ), \\
g_n(\theta,r)&=(\theta+[n q_n^{\si}] r, r), \\
h_n&= g_n \ci \phi_n, \quad H_n =h_1 \circ \dots \circ  h_n, \\
f_n&=H_n \ci R_{\a_{n+1}} \ci H_n^{-1}.
\end{aligned}
\end{equation}

\subsection{Proof of convergence.}
Let $\a, \delta$ and $\sigma$ be as in the statement of Theorem \ref{analytiques}, and let $\rho>0$ be 
fixed. Let $\a_n=p_n/q_n$ be a sequence such that $|\a-\a_n| $ is decreasing and 

\begin{itemize}
\item[(P1)]  For all $n \in \NN$, 
$$
|\a-\a_n|<\exp (-q_n^{1+3\si}).
$$

\end{itemize}
By eventually extracting from $\a_n$ we can assume that this sequence  also  has the following properties:

\begin{itemize}
\item[(P2)]
Denote the  lift of the inverse of the diffeomorphism 
$H_{n}$ from (\ref{eqdef_apprintro})
by $((H_n^{-1})_1,\,(H_n^{-1})_2)$, and set 
$$
\rho_n:=\max_{i=1,2}  \inf_{p\in\ZZ}
\|(H_n^{-1})_i+p\|_\rho, \quad 
\rho_0: =\rho.
$$
Then for  all $n \in \NN$,
$$
q_{n}^\si\geq 4\pi n \rho_{n-1}+ \ln(8\pi n q_n^{\si + 4}).
$$

\item[(P3)] With the definition of $\|DH\|_\rho$ of Section \ref{analtop}, 
we have for all $n \in \NN$, and for all $t$ such that $|t-\a|\leq |\a_n-\a|$,
$$
q_n\geq \|D(H_{n-1})R_{t} \circ H_{n-1}^{-1}\|_\rho.
$$

\item[(P4)] For all $n \in \NN$ 
$$
\|D(H_{n-1})\|_0 \leq \ln q_n.
$$

\end{itemize}

Properties (P2)--(P4) are  possible to guarantee by choosing $q_n$ sufficiently large because
$H_{n-1}$ does not depend on $q_n$.

The first three properties are used to prove the convergence, and the latter one is estimate 
(\ref{estimate}), needed for the
proof of weak mixing of the limit diffeomorphism, which will be done with the help of Proposition 
\ref{criterion}.

The following statement implies the convergence of the 
sequence $f_n$.

\begin{lemm}\label{lco} Suppose  $\a_n=\frac{p_n}{q_n}$
satisfies  (P1)--(P3) 
for some fixed $\si >0$ and $\rho>0$. 
Then, for any $n$ large enough, we have:

\noindent (a) the diffeomorphisms defined by (\ref{eqdef_apprintro}) satisfy:
$$
d_\rho(f_n,f_{n-1}) \leq \exp (-q_{n});
$$
(b) for any $m\leq q_{n+1}$ it holds:
$$
d_0(f_n^m,f^m)\leq  \frac{1}{2^n}.
$$
\end{lemm}

\begin{proof}
With  the notations above, using the Mean value theorem and (P3), we
have (for some $t$ between $\a_{n}$ and $\a_{n+1}$):
\begin{equation}\label{eqec}
\begin{aligned}
d_\rho(f_n,f_{n-1}) &\leq \|(D H_{n-1})R_{t} \ci H_{n-1}^{-1}\|_\rho 
\|(h_n\ci R_{\a_{n+1}} \ci h_n^{-1} - R_{\a_n}) \ci  H_{n-1}^{-1}\|_{\rho} \\
& \leq q_n
\|h_n\ci R_{\a_{n+1}} \ci h_n^{-1} -R_{\a_n}\|_{\rho_{n-1}}.
\end{aligned}
\end{equation}

Denote $(\cos2\pi q_n (z+\a_{n+1})-\cos2\pi q_n z)$ by $R(z)$. 
For an arbitrary $s \geq 0$, we can write:
\begin{equation}\label{eqR}
\begin{aligned}
\|R\|_s  \leq & \|e^{2\pi i q_n z} \|_{s} |1-e^{2\pi i q_n \a_{n+1}}| \leq 
2 \pi q_n \| e^{2\pi i q_n z} \|_{s} |\a_{n+1}-\a_n|  \\
\leq 4 & \pi q_n \| e^{2\pi i q_n z} \|_{s} |\a-\a_n|,
\end{aligned}
\end{equation}
(we used the estimate $|\a_{n+1}-\a_n| \leq 2 |\a-\a_n|$).
By the definition of $h_n$, 
$$
h_n \ci R_{\a_{n+1}} \ci  h_n^{-1} - R_{\a_n} =
([nq_n^{\si}] q_n^2 R(\theta -[nq_n^{\si}] r)+ (\a_{n+1}-\a_n)\, ,\, q_n^2 R(\theta-[nq_n^{\si}] r)).
$$
Then
$$
\|h_n \ci R_{\a_{n+1}} \ci  h_n^{-1} - R_{\a_n}\|_s 
\leq 2nq_n^{2+\si}\|R(\theta -[nq_n^{\si}] r)\|_s
.
$$
By (\ref{eqR}), it is less than 
\begin{equation}\label{otsenka1}
 8\pi n q_n^{3+\si} 
\|\exp(2\pi i q_n (\theta -[nq_n^{\si}] r) )\|_s  
 |\a-\a_n| .
\end{equation}
Applying  (\ref{eqec}), (\ref{otsenka1}), (P2) and (P1) in sequence, we get: 
$$
\begin{aligned}
d_\rho(f_n,f_{n-1}) \leq &
q_n \|h_n\ci R_{\a_{n+1}} \ci h_n^{-1} -R_{\a_n}\|_{\rho_{n-1}}  \\
\leq & 8 \pi n q_n^{4+\si} \exp(4\pi nq_n^{1+\si} \rho_{n-1} ) |\a-\a_n|\leq 
\exp(q_n^{1+2\si} ) |\a-\a_n| \\
\leq & \exp(q_n^{1+2\si}(1 -q_n^{\si}))
 \leq \exp(-q_n^{1+2\si})<\exp(-q_n).
\end{aligned}
$$

The second part of the claim is proved in the same way. One has to
note that $f_n^m=h_n \ci S^m_{\a_{n+1}} \ci  h_n^{-1}=h_n \ci
R_{m\a_{n+1}}\ci  h_n^{-1} $, and 
$$
d_0(f^m,f_n^m)=\sum\limits_{j=n}^\infty d_0(f_j^m,f_{j+1}^m).
$$
\end{proof}


\subsection{Proof of weak mixing}

For the proof of weak mixing, we shall use Proposition  \ref{criterion} that was
proved in the previous section. In order to apply the lemma, we choose
a sequence ($m_n$), $m_n\leq q_{n+1}$ (in this case, by Lemma 
\ref{lco} (b), (\ref{zvezdochka}) holds), and a sequence of standard partial
decompositions ($\eta_n$) consisting of horizontal intervals with length less than $1 / q_n$, $\eta_n \to 
\epsilon$, such that the
diffeomorphism
\begin{equation}\label{oprPhi}
\Phi_n:=\phi_n \circ R_{\a_{n+1}}^{m_n}\circ \phi_n^{-1}
\end{equation}
$(\frac{1}{nq_n}, 0, \frac{1}{n} )$-distributes any 
interval  $I_n\in\eta_n$. 

\subsubsection{Choice of the mixing sequence $m_n$.}
We shall assume that 
$$
q_{n+1}\geq q_n^7.
$$
Define 
$$
m_n=\min \left\{ m\leq q_{n+1} \mid 
\inf_{k\in \ZZ} 
\left| m\frac{q_{n}p_{n+1}}{q_{n+1}}-1/2 + k \right| < 
\frac{q_n}{q_{n+1}}   \right\}. 
$$
Note that the set of numbers $m$ above is non-empty. Indeed,
since $p_{n+1}$ and $q_{n+1}$ are relatively prime, the set 
$\{j\,\frac{q_{n}p_{n+1}}{q_{n+1}}\mid j=0,\dots q_{n+1} \}$ on the circle
contains $\frac{q_{n+1}}{GCD(q_n,q_{n+1})}$, which is at least 
$\frac{q_{n+1}}{q_n}$, different equally distributed points.

We shall use the following estimate, which follows from the above assumption
on the growth of $q_n$: 
\begin{equation}\label{analm}
|m_n q_n \a_{n+1} -1/2|(\text{mod}\,1)  
\leq \frac{q_n}{q_{n+1}}\leq q_n^{-6}.
\end{equation}


\subsubsection{Stretching of the diffeomorphisms $\Phi_n$.}\label{diffPhi}
Consider the set 
\begin{equation}\label{setB_n} 
B_n=\bigcup_{k=0}^{2q_n} 
\left[\frac{k}{2q_n}-\frac{1}{2q_n^{3/2}}\,,\,\frac{k}{2q_n}+\frac{1}{2q_n^{3/2}}\right].
\end{equation}
We shall see that $\Phi_n$ displays strong stretching in the vertical
direction on small
horizontal intervals, lying outside $B_n$. To do this, we shall use
the notion  of
uniform stretch from \cite{F}, which we recall here.  

\begin{defi}[Uniform stretch]\label{dus}
Given $\e>0$ and $k>0$, we say that a real function $f$ on an interval 
$I$ is ($\e,k$)-uniformly stretching on $I$ if for $J=
[\inf_{I}f,\sup_{I} f]$ 
$$
\la(J) \geq k,
$$
and for any interval $\tilde J \subset  J$ we have:
$$
\left|\frac{ \la(I\cap f^{-1} \tilde J)}{\la(I)}-
\frac{\la(\tilde J)}{ \la(J)}\right| 
\leq \e \frac{ \la(\tilde J)}{\la( J)}.
$$
\end{defi}

The following criterion,  that is easy to verify, 
gives a necessary and
sufficient condition for a real function (of class at least $C^2$) 
to be uniformly stretching. The proof can be found in \cite{F}. 
\begin{lemm}[Criterion for uniform stretch]\label{lus}
If $f$ satisfies:
\begin{eqnarray*}
\inf_{x\in I}|f'(x)| \la(I) &\geq& k,  \\
  \sup_{x\in I}|f''(x)|  \la(I) &\leq& \e \inf_{I}|f'(x)|, 
\end{eqnarray*}

then $f$ is ($\e,k$)-uniformly stretching on $I$.
\end{lemm}


\begin{lemm}\label{cm}
Under the conditions of Theorem \ref{analytiques}, the transformation  
$\Phi_n$  
has a lift of the form:
$$ 
\Phi_n(\theta,r) =(\theta+m_n\a_{n+1}, \ r+\psi_n(\theta)),
$$
where $\psi_n$ satisfies: 
\begin{equation}\label{proizv}
\inf_{\TT \setminus B_n}|\psi_n'| \geq q_n^{5/2},
\quad \sup_{\TT \setminus B_n}|\psi_n''| \leq 9 \pi^2 {q_n}^4.
\end{equation}
\end{lemm}

\begin{proof}
By definition, 
$\Phi_n$ 
has the desired form with 
$$
\psi_n=q_n^2 (\cos (2\pi (q_n \theta +m_nq_n\a_{n+1})) - \cos (2\pi q_n
\theta) ) = -2 q_n^2 \cos (2\pi q_n \theta)+\si_n,
$$ where
$$
\si_n=q_n^2(\cos( 2\pi (q_n\theta +m_nq_n\a_{n+1})) -
\cos (2\pi (q_n\theta +1/2))). 
$$

With the help of the Mean value theorem and estimate
(\ref{analm}), one easily verifies that
$|\si_n'|<1$, and 
$|\si_n''|<1$. 

Note that $B_n$ are chosen in such a way that 
$$
\inf_{\TT\setminus B_n}
|\sin(2\pi q_n \theta ) | \geq q_n^{-1/2}.
$$
The statement follows by calculation.
\end{proof}


\subsubsection{Choice of the decompositions $\eta_n$.}
Let us define a standard
partial decompositions $\eta_n$ of $\TT^2$, meeting the
conditions of Proposition \ref{criterion}.

Let $\hat \eta_n=\{I_n\}$ be the partial decomposition of 
$\TT \setminus B_n$, containing all the intervals  $I_n$ such that
$$
\psi_n(I_n)=[0,1) \mod 1.
$$
We define $\eta_n=\{I\times \{r\}\mid I \in \hat\eta_n,\  r\in
\TT \}$. Note that, for any $I_n \in \eta_n$, we have: $\pi_r(\Phi(I_n))=\TT$.
\begin{lemm}
Let  $\eta_n$ be defined as above.
Then,  for any $I_n \in \eta_n$, 
$$
\la(I_n)\leq q_n^{-5/2},
$$  
and $\eta_n\to\epsilon$.
\end{lemm}
\begin{proof}

By Lemma \ref{cm}, 
$\inf_{\TT \setminus B_n}|\psi_n'| \geq q_n^{5/2} $.
Therefore,  
$\la(I_n)\leq q_n^{-5/2}$ for any $I_n\in \eta_n$.

Since the diameter of the atoms of $\eta_n$ goes to zero when $n$ grows, it is enough to show that
the total measure of the decompositions goes to 1 when $n$ grows.
The total measure of $\eta_n$ equals: 
$$
\begin{aligned}
\sum\limits_{I_n\in \hat \eta_n}\la (I_n) 
\leq 1-\la(B_n)-4q_n\max_{I_n\in \hat \eta_n}\la (I_n) \\ 
\leq 1 - 2 q_n (q_n^{-3/2} + 2q_n^{-5/2})< 1-3q_n^{-1/2}\to 1.
\end{aligned}
$$
\end{proof}


\subsubsection{Proof of weak mixing.}
To prove weak mixing of $f$, we shall apply Proposition
\ref{criterion}. Since (\ref{zvezdochka}) holds by Lemma \ref{lco},
estimate (\ref{estimate}) holds by Property (P4), the sequence of decompositions
$\eta_n\to \epsilon$ by the lemma above, it is left to  
verify condition (2) of Proposition \ref{criterion}, which we pass to.
\begin{lemm}
Let $I_n\in\eta_n$, $\Phi_n$ be as in (\ref{oprPhi}). 
Then $\Phi_n(I_n)$ is  $(\frac{1}{nq_n}, 0, \frac{1}{n} )$-distributed.
\end{lemm}

\begin{proof}
By the choice of $\eta_n$, $\pi_r( \Phi_n(I_n))=\TT$, and hence, $\de$ in the definition of
($\ga,\de,\e$)-distribution can be taken equal to 0. 

We have seen that $\Phi_n$ has a lift  $\Phi_n(r,\theta)=(\theta +m_n\a, r+
\psi_n(\theta))$. Hence, $\Phi_n(I_n)$ is contained in the
vertical strip $(I_n+m_n\a)\times\TT$.
By the lemma above, $\la(I_n)\leq \frac{1}{q_n^{5/2}}<\frac{1}{nq_n}$ for any $I_n\in \eta_n$.
Hence, we can take $\ga =\frac{1}{nq_n}$.

Our fixed $I_n$ has the form $I\times \{r\}$ for some
$r\in\TT$ and $I \in \hat \eta_n$.
For any $J\subset \TT$, the fact that
$
\Phi_n (\theta, r)\in \TT \times J 
$
is equivalent to 
$
\psi_n(\theta)\in J-r 
$.
Lemma \ref{cm} implies the estimate:
$$
\frac{\sup_{I_n\in\eta_n} |\psi_n''| }{\inf_{I_n\in\eta_n}
  |\psi_n'|}\la(I_n)\leq \frac{9 \pi^2}{q_n}<\frac{1}{n}.
$$
Then, by Lemma \ref{lus} (Criterion for uniform stretch), $\psi_n$ is
($\frac{1}{n},1$)-uniformly stretching. Hence,
for any interval $J\subset \TT$, 
the following holds:
$$
\begin{aligned}
|\la (I_n\cap \Phi_n^{-1} (\TT\times J))-\la (I_n)\la ( J)|
=&|\la (I \cap \psi_n^{-1} (J-r))-\la (I_n)\la ( J)| \\
\leq & \frac{1}{n}
\lambda(I_n)\la(J),
\end{aligned}
$$
and we take $\e=\frac{1}{n}$ in the definition of
($\ga,\de,\e$)-distribution. 
\end{proof}

We have shown that $\Phi_n$ and $\eta_n$ verify the conditions of Proposition
\ref{criterion}.  It implies that $f$ is weak mixing.

\section{$C^\infty$-case on the torus, annulus and disc}
\label{section-smooth}

Sections 5.1--5.4 are devoted to $M=\AAA$ and $M=\TT^2$. The case of
the disc $\DD^2$ is studied in Section 5.5.

\subsection{Statement of the result.} 
Take any $0<\si < 1$. 
On $M=\AAA$, consider the following transformations:
\begin{equation} \label{eqsmooth_appr}
\begin{aligned} 
&g_n(x,y)=(x+[nq_n^\sigma] y,\ y), \\
&h_n= g_n \ci \phi_n, \quad H_n =h_1 \ci \dots \ci  h_n,\\
&f_n=H_n \ci R_{\a_{n+1}} \ci H_n^{-1};
\end{aligned}
\end{equation}
where  the sequence $\a_n= p_n/q_n$, converging to $\a$, and 
 the diffeomorphisms  $\phi_n$, satisfying 
\begin{equation}\label{com5}
R_{\frac{1}{q_n}} \circ \phi_n = \phi_n\circ R_{\frac{1}{q_n}},
\end{equation}
will be constructed  in Section
\ref{sectionphi} below so that
\begin{theo}\label{thsmooth} 
For any Liouville number $\a$, there exists a sequence 
$\a_n$ of rationals and a sequence $\phi_n$ of measure preserving diffeomorphisms 
satisfying (\ref{com5}) such that the diffeomorphisms $f_n$,
constructed as in (\ref{eqsmooth_appr}), converge in the sense of  the  $\diff^\infty(M)$ topology, 
the limit diffeomorphism 
$f=\lim\limits_{n\to \infty} f_n$ being weak mixing and $f\in \A_\a(M)$.
Moreover, for any $\e>0$, the parameters can be chosen 
so that 
$$
d_\infty(f , R_{\a})<\e.
$$
\end{theo} 
 
\begin{rema}\label{posl}
This result implies Theorem \ref{liouville}. Indeed, it follows directly from  
Theorem \ref{thsmooth}, that
weak mixing diffeomorphisms are dense in $\A_\a(M)$. It is a general
fact (see \cite{H}) that, in this case, weak mixing diffeomorphisms are generic in
$\A_\a(M)$ with our topology.
\end{rema}

\subsection{Construction of $\phi_n$.}\label{sectionphi}
We begin by constructing a ``standard diffeomorphism" on the square 
$[-1,1]\times [-1,1]=[-1,1]^2$, from which $\phi_n$ will be obtained
by a rescaling of the domain of definition.
\subsubsection{Preliminary construction.}
For a fixed $\e<1/2$, 
consider the squares $\Delta=[-1,1]^2$,  
$\Delta(\e)= [-1+\e,1-\e]^2$ and $\Delta(2\e)$. 
\begin{lemm}\label{lsmooth1}
For any $\e<1/2$ 
there exists a smooth measure-preserving diffeomorphism
$\ph=\ph (\e)$ of $\mathbb R^2$, equal to the identity outside 
$\Delta (\e)$ 
and rotating the square $\Delta(2\e)$
by $\pi/2$.
\end{lemm}

\begin{proof}
Let $\psi=\psi(\e)$ be a smooth transformation satisfying 
$$
\psi(\theta,r)=
\begin{cases}
(\theta,r) &\text{ on } \mathbb R^2 - \Delta(\e), \\
(\theta/5, r/5) &\text{ on } \Delta(2\e),
\end{cases}
$$
and $\eta$ be a smooth transformation, such that
$$
\eta(\theta,r) =
\begin{cases}
(r,-\theta) &\text{ on }  \{\theta^2+ r^2\leq 1/3\}, \\
(\theta,r) &\text{ on } \{\theta^2+ r^2\geq 2/3\}.
\end{cases}
$$
Then the composition 
$$
\tilde \ph:=\psi^{-1} \eta \psi
$$
provides the desired geometry. 
Moreover, it preserves the 
Lebesgue measure on the set  
$$
U=(\mathbb R^2 -\Delta(\e))\cup \Delta(2\e).
$$
However, it does not have to preserve the area on the whole of $\Delta$. We describe now a deformation 
argument following Moser \cite{moser} that provides 
an {\it area-preserving} diffeomorphism $\ph$ on $\Delta$,
coinciding with $\tilde \ph$ on $U$.

Let $\Omega_0$ denote the usual volume form on $\mathbb R^2$, 
and consider $\Omega_1:=\tilde\ph^* \Omega_0$. 
We shall find a diffeomorphism 
$\nu$ equal to the identity on the set $U$,
and such that $\nu^* \Omega_1=\Omega_0$. 

Let $\Omega'=\Omega_1-\Omega_0$, and note that 
$\Omega'=d(\omega_0-\tilde\ph^* \omega_0)$, where
$\omega_0$ is the standard 1-form $\frac{1}{2}(\theta dr-r d\theta)$.
Consider the volume form 
$$
\Omega_t=\Omega_0+t\Omega'.
$$
Since it is non-degenerate, there exists a unique vector
field $X_t$ such that 
\begin{equation}\label{eqforms}
\Omega_t(X_t,\cdot)=(\omega_0-\tilde \ph^* \omega_0)(\cdot).
\end{equation}

One can integrate the obtained vector field to get the 
one-parameter family of diffeomorphisms $\{\nu_t\}_{t\in [0,1]}$,
$\dot \nu_t=X_t(\nu_t)$, $\nu_0=id$. Then $\nu=\nu_1$ is the desired
coordinate change. Indeed, one verifies by calculation that
$$
\frac{d}{dt}\nu_t^*\Omega_t=0.
$$
Hence, $\nu_1^*\Omega_1=\nu_0^*\Omega_0=\Omega_0$.

By an explicit verification, one obtained 
that $\tilde\ph^*$ preserves the form $\omega_0$
on $U$ (for this note that $\tilde\ph$ on $U$ is an explicit linear
transformation). Then on $U$ equation (\ref{eqforms}) writes as
$\Omega_t(X_t,\cdot)=0$.
Since $\Omega_t$ is non-degenerate,
this implies that $X_t=0$ on $U$, hence $\nu=\nu_0=id$ on $U$, as claimed.
The desired area-preserving diffeomorphism is 
$$
\ph=\nu \tilde \ph.
$$
\end{proof}

\subsubsection{Construction of $\phi_n$} \label{consphi}

Let us first define 
$\phi_n$ on the fundamental domain 
$D_n=[0,1/q_n]\times[0,1] $. The line $\theta=1/2q_n$ divides $D_n$ into halves: $D_n^1=[0,1/(2q_n)] \times [0,1]$ and
$D_n^2=(1/(2q_n),1/q_n) \times [0,1]$. 
On $D_n^1$, consider the affine transformation
$C_n(\theta,r)=(4q_n \theta-1,2r-1)$, sending $D_n^1$ onto the square
$\Delta=[-1,1]^2$. Let $\ph_n$ be the diffeomorphism given by Lemma \ref{lsmooth1}
with $\e=1/(3n)$,
and set
\begin{equation}\label{constrphi}
\phi_n:=C_n^{-1}\circ \ph_n \circ  C_n. 
\end{equation}
We define $\phi_n=Id$ on $D_n^2$.
Note that $\phi_n$ is smooth and area-preserving on $D_n$, and equals identity on the
boundary of $D$. 
We extend it periodically to the whole $\RR^2$ by the formula:
$$
\phi_n \circ R_{\frac{1}{q_n}}= R_{\frac{1}{q_n}}\circ \phi_n,
\quad \phi_n(\theta,r+1)=\phi_n(\theta,r) + (0,1).
$$ 
The transformation 
$\phi_n$, defined in this way, becomes a diffeomorphism both on $\TT^2$ and on
$\AAA$ in a natural way.
  
For a fixed $n$, let us denote by $D_{n,j}$ and $D_{n,j}^i$ 
(for $i=1,2$, $j \in \ZZ$) the
shifts of the fundamental domain $D_n$ of $\phi_n$:
$$
D_{n,j+q_n}=D_{n,j}=R_{\frac{j}{q_n}} (D_n), \text{ and } D_{n,j+q_n}^i= D_{n,j}^i=R_{\frac{j}{q_n}} 
(D_n^i).
$$

\subsubsection{Notation}\label{notation}
For a diffeomorphism $F$ of $M$ (not necessarily homotopic to the
identity), we shall denote by the  
same letter its lift of the form:
$$
F(x,y)=(ax+by+f_1(x,y),cx+dy+f_2(x,y)),
$$
where $f_i:\RR^2\to \RR$ are, in the case of the torus,
$\ZZ^2$-periodic with the property
$ \|f_i\|_{0}= \inf_{p\in \ZZ}\|f_i+p\|_{0}$; and for the case of
the annulus, $f_i$ are $\ZZ$-periodic in the first component, and such that
$ \|f_1\|_{0}= \inf_{p\in \ZZ}\|f_1+p\|_{0}$. Note that the
diffeomorphisms in our constructions are defined by their lifts,
satisfying this property.
For $k$-smooth diffeomorphisms $F: \RR^2 \rightarrow \RR^2$ 
we define by $F_i$ the $i$-th coordinate
function, and denote
$$
|\!|\!| F |\!|\!|_{k} := 
\max \{ \| D_a F_i \|_{0},\, \| D_a (F^{-1})_i
\|_{C^0} \mid i=1,2,\ \  0 \leq
|a| \leq k \}.  
$$
 
\subsubsection{Discussion of the properties of $\phi_n$} \label{ppphi}

\begin{figure}
\psfrag{d}[][][0.7]{$D_{n,j}^1$}
\psfrag{dd}[][][0.7]{$D_{n,j}^2$}
\psfrag{j}[][][0.7]{$I$}
\psfrag{phi}[][][0.7]{$\phi_n(I)$}
\psfrag{phi=id}[][][0.7]{$\phi_n=\text{Id}$}
\includegraphics[width=3.5cm]{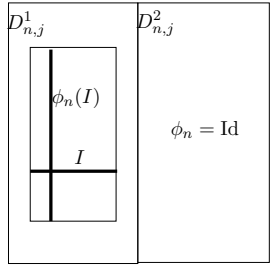}
\caption{Action of $\phi_n$}\label{kartina_phi}
\end{figure}

We have constructed $\phi_n$ so that $\phi_n$ equals identity on $D_{n,j}^2$, $j \in \ZZ$, and on 
$D_{n,j}^1$ the image of any 
interval ${I}_{n,j} \times\{r\}$, where
$r\in[1/(3n),1-1/(3n)]$, and
\begin{equation}\label{eqint_smooth}
{I}_{n,j} = \left[\frac{j}{q_n}+\frac{1}{6nq_n},
\frac{j}{q_n} + \frac{1}{2q_n}-\frac{1}{6nq_n}\right], 
\end{equation} 
with $j=0,\dots q_n-1$, 
both under $\phi_n$ and $\phi_n^{-1}$, is an interval 
of type $\{\theta\} \times[1/(3n),1-1/(3n)]$
for some $\theta \in I_{n,j}$ (see Figure \ref{kartina_phi}).

Moreover, the following holds:
\begin{lemm}\label{lsmooth2}
For all $k\in\NN$ the diffeomorphisms $\phi_n$ constructed above 
satisfy: 
$$
|\!|\!| \phi_n |\!|\!|_{k} \leq c(n,k) q_n^{k},
$$
where $c(n,k)$ is independent of $q_n$.
\end{lemm}  
\begin{proof}
The desired estimate follows from (\ref{constrphi}) by the product
rule (it is important that $\ph_n$ is independent of $q_n$). 
\end{proof}
\begin{rema}\label{remphi} For any $n$, the  
construction implies that 
$\phi_n(\theta,r)=Id \, $ in the domains
$0\leq r < 1/(6n)$ and $1-1/(6n)< r \leq 1$. It is easy to
verify that in the same domains diffeomorphisms $f_n$ from 
(\ref{eqsmooth_appr}) equal $R_{\a_{n+1}}$. 
\end{rema}


\subsection{Proof of convergence}
In the proof we shall use the following lemma:
\begin{lemm}\label{Alistair}
Let $k \in \NN$, and 
$h$ be a diffeomorphism of $M$. Then for all $\alpha, \beta
\in \RR$ we obtain
\begin{align}\label{eq:MainLem}
d_k( h R_{\alpha} h^{-1}& , h R_{\beta} h^{-1})
\leq C_k   |\!|\!|h|\!|\!|_{k+1}^{k+1} | \alpha - \beta|,
\end{align}
where $C_k$ only depends on $k$, and $C_0=1$.
\end{lemm}

\begin{proof}
We give the proof for the case $M=\TT^2$; for the annulus, the proof is
obtained by minor modifications.
Note that $D_a h_i$ for $|a|\geq 1$ is
$\ZZ^2$-periodic. Hence, for any $g:\RR^2\to \RR^2$, we have: 
$\sup_{0< x,y < 1} |(D_ah_i)(g(x,y))|\leq 
|\!|\!|h|\!|\!|_{|a|}$. 

For $k=0$, the statement 
of the lemma follows directly from the  Mean value theorem. 

We claim that for $j$ with $|j|=k$ the partial derivative 
$D_{j} ( h_i  R_{\alpha} h^{-1}  - h_i R_{\beta} h^{-1} )$ will
consist of a sum of terms with each term being the product of a
single partial derivative
\begin{equation}\label{eq:term1}
\bigl( D_a h_i \bigr) (R_{\alpha} h^{-1} ) -\bigl( D_a h_i\bigr)
    (R_{\beta} h^{-1})
\end{equation}
with $|a| \leq k$, and at most $k$ partial derivatives of the form 
\begin{equation}\label{eq:term2}
D_b h^{-1}_j 
\end{equation}
with $|b| \leq k$. This clearly holds for $k=1$. We proceed by
induction. 

By the product rule we need only consider the effect of
differentiating (\ref{eq:term1}) and (\ref{eq:term2}).
Applying $D_c$ with $|c|=1$ to (\ref{eq:term1}) we get:
$$
\sum_{|b|=1} \Bigl(\bigl( D_b D_a h_i \bigr) (R_{\alpha} h^{-1} )
- \bigl( D_b D_a h_i \bigr) (R_{\beta} h^{-1}) \Bigr) D_c h^{-1}_b,
$$
which increases the number of terms of the form (\ref{eq:term2}) in the
product by
1.  Differentiating (\ref{eq:term2}) we get another term of the form
(\ref{eq:term2}) but with $|b| \leq k+1$. 
  
Now we  estimate:
\begin{gather*}
\| \bigl( D_a h_i \bigr) (R_{\alpha} h^{-1} )- \bigl( D_a h_i
\bigr) (R_{\beta}
h^{-1} )\|_0 \leq |\!|\!|  h |\!|\!|_{|a|+1} | \alpha - \beta |, \\
\| D_c h^{-1}_j \|_0 \leq |\!|\!| h |\!|\!|_{|c|}.
\end{gather*}
Taking the
inverse maps and applying the result we just proved gives
\eqref{eq:MainLem}.
\end{proof}

\begin{lemm}\label{convergence} 
For an arbitrary $\e>0$, 
let $k_n$ be a growing sequence of natural numbers, such that
$\sum_{n=1}^\infty 1/k_n<\e$. 
Suppose that, in construction (\ref{eqsmooth_appr}), 
we have: $|\a-\a_1|<\e$ and for any $n$ 
\begin{eqnarray} \label{conv}  
|\a - \a_{n}| <  \frac{1}{ 2 k_n C_{k_n}{|\!|\!|H_n|\!|\!|}_{k_n+1}^{k_n+1}},
\end{eqnarray}
where $C_{k_n}$ are the constants from Lemma \ref{Alistair}.
Then the diffeomorphisms 
$f_n = H_n\circ R_{\a_{n+1}}\circ H_n^{-1}$ converge in the
$\diff^\infty$ 
topology to a measure preserving diffeomorphism $f$, and
$$
d_\infty(f , R_{\a})<3\e.
$$ 

Moreover, the sequence of diffeomorphisms
\begin{eqnarray} \label{ftilde}
\hat{f}_n:=H_{n}\circ R_{\a }\circ H_{n}^{-1} \in \A_\a
\end{eqnarray}
also converges to $f$ in the $\diff^\infty$ 
topology, hence $f\in \A_\a$.

\noindent Furthermore, if for a sequence of positive integers $m_n$ we have for
all $n$: 
\begin{eqnarray} \label{conv2}  
|\a - \a_{n}| <  \frac{1}{2^{n+1} m_{n-1} |\!|\!|H_n|\!|\!|_{1}  },
\end{eqnarray}
then for any $m \leq m_n$ we have 
\begin{eqnarray} \label{conv3} 
d_0(f^m , f^m_n) \leq \frac{1}{2^n}. 
\end{eqnarray} 
\end{lemm}

\begin{proof}

By construction we have: $h_n \circ R_{\a_n} = R_{\a_n} \circ h_n$. Hence,
$$
f_{n-1}=H_{n-1}\circ R_{\a_{n}}\circ H_{n-1}^{-1} =
H_{n}\circ R_{\a_{n}}\circ H_{n}^{-1}.
$$
By Lemma \ref{Alistair}, for all $k$ and $n$, 
$$
\begin{aligned}
d_k(f_n , f_{n-1})=&d_k(H_n\circ R_{\a_{n+1}}\circ H_n^{-1} ,
H_{n}\circ R_{\a_{n}}\circ H_{n}^{-1}) \\ 
\leq &C_k |\!|\!| H_{n} |\!|\!|_{k+1}^{k+1} 
| \a_{n+1} - \a_n|. 
\end{aligned}
$$ 
Estimating 
$| \a_{n+1} - \a_n|\leq 2 |\a - \a_n| $, and using assumption
(\ref{conv}), we get for any $k\leq k_n$: 
$$
d_k(f_n , f_{n-1}) \leq d_{k_n}(f_n , f_{n-1}) \leq
\frac{2 C_{k_n} |\!|\!| H_n |\!|\!|_{{k_n}+1}^{{k_n}+1}} 
{ 2{k_n}C_{k_n} |\!|\!| H_n|\!|\!|_{{k_n}+1}^{{k_n}+1} }
\leq \frac{1}{k_n}.
$$
Hence, for any fixed $k$, the sequence ($f_n$) converges in
$\diff^k$, and therefore, in $\diff^\infty$. 
Moreover, one easily computes (using the definition of the
$d_\infty$-metric) that 
$$
d_\infty(f , R_\a)\leq |\a-\a_1|+
\sum_{n=1}^\infty d_\infty(f_n , f_{n-1})
<3\e;
$$
(here we denoted $f_0=R_{\a_1}$).

To prove that $f\in \A_\a$, we show that the sequence of functions 
$\hat{f}_n\in \A_\a$ converge to $f$. For this it is enough to
note that, for any $n$ and $k\leq k_n$, Lemma \ref{Alistair} and assumption
(\ref{conv}) imply:
$$
\begin{aligned}
d_k(f_n,\hat{f}_n)=&d_k(H_n\circ R_{\a_{n+1}}\circ H_n^{-1} ,
H_{n}\circ R_{\a }\circ H_{n}^{-1}) \\
\leq &C_{k_n} |\!|\!| H_{n} |\!|\!|_{k_n+1}^{k_n+1} 
| \a_{n+1} - \a|\leq \frac{1}{k_n}.
\end{aligned}
$$

To prove the third statement of the lemma, note that
for any $m\leq m_{n-1}$,
$$  
\begin{aligned}
d_0(f_n^m,f_{n-1}^m) =& d_0( H_n\circ  R_{m\a_{n+1}}\circ  H_n^{-1},
H_{n}\circ R_{m\a_{n}}\circ H_{n}^{-1} )  \\
\leq &|\!|\!| H_n|\!|\!|_{1} 2m |\a - \a_n| \leq \frac{1}{2^{n}}.
\end{aligned}
$$
Then $d_0( f^m,f_{n-1}^m ) \leq 
\sum\limits_{i=n}^\infty d_0( f_i^m, f_{i-1}^m) = \frac{1}{2^{n-1}}$.
\end{proof}

Let a Liouville number $\a$ be fixed. Here we show that, for any given
sequence $k_n$, the sequence
of convergents $\a_n$ of $\a$ can be chosen so that (\ref{conv})
holds, and for any $m_{n-1}\leq q_{n}$, (\ref{conv2}) holds.

\begin{lemm}\label{Nnda}
Fix an increasing sequence $k_n$ of natural numbers, satisfying
$\sum_{n=1}^\infty 1/k_n<\infty$, and let the constants $C_n$ be as in
Lemma \ref{Alistair}. 
For any Liouville number $\a$, there exists a sequence of convergents
$\a_n= p_n/q_n$, such that the diffeomorphisms 
$H_n$,
constructed as in (\ref{eqsmooth_appr}) 
with these $\a_n$ and with $\phi_n$ given by (\ref{constrphi}),
satisfy (\ref{conv}) and (\ref{conv2}) 
with any $m_{n-1}\leq q_{n}$. Further, we can choose $\a_n$ so that in addition (\ref{estimate})) holds.
\end{lemm}
\begin{proof}
By Lemma \ref{lsmooth2}, we have: 
$|\!|\!|\phi_n|\!|\!|_{k}\leq c_1(n,k) q_n^k$.
Then for $h_n$ as in (\ref{eqsmooth_appr}), we get:
$$
|\!|\!| h_{n} |\!|\!|_{k}\leq c_2(n,k) q_n^{2k}.
$$
With the help of the Faa di Bruno's formula 
(that gives an explicit equation for the $n$-th derivative
of the composition), we estimate: 
$$
|\!|\!| H_{n} |\!|\!|_{k}\leq 
|\!|\!| H_{n-1}\ci h_n |\!|\!|_{k}\leq c_3(n,k)q_n^{2k^2},
$$
where $c_3(n,k)$ depends on the derivatives of $H_{n-1}$ up to
order $k$, which do not depend on $q_n$.
Suppose that, for each $n$, $q_n$ is chosen so that 
$$
q_n\geq c_3(n,n+1).
$$
Then $|\!|\!|H_{n}|\!|\!|_{{k_n}+1}\leq q_n^{2({k_n}+1)^2+1}
\leq q_n^{3({k_n}+1)^2}$.
We choose the sequence of convergents of $\a$ satisfying
$$
|\a-\a_n|=|\a-p_n/q_n| <\frac{1}{2^{n+1} k_n C_{k_n} q_n^{3({k_n}+1)^3+1}};
$$
the latter is possible since $\a$ is Liouville.
Then 
$$
|\a-\a_n|<\frac{1}{2^{n+1} q_n k_n C_{k_n} 
|\!|\!|H_{n}|\!|\!|_{{k_n}+1}^{{k_n}+1}},
$$ 
which implies both (\ref{conv}) and (\ref{conv2}). As for (\ref{estimate}), i.e.  $\|DH_{n-1}\|_0\leq \ln 
q_n$, it is possible to have it just by choosing $q_n$ large enough.
\end{proof}

\subsection{Proof of weak mixing}

\subsubsection{Choice of the mixing sequence $m_n$.}
We shall assume that for all $n$ we have:
\begin{equation}\label{grocondsmooth}
q_{n+1}\geq 10 n^2q_n.
\end{equation}
Define, as in the analytic case, 
$$
m_n=\min \left\{ m\leq q_{n+1} \mid \inf_{k\in \ZZ}
\left| m\frac{q_{n}p_{n+1}}{q_{n+1}}-1/2+ k \right|\leq 
\frac{q_n}{q_{n+1}}   \right\}. 
$$

Let $a_n=(m_n\a_{n+1}-\frac{1}{2q_n})\mod \frac{1}{q_n}$. Then the choice of $m_n$ and
the growth condition (\ref{grocondsmooth}) imply:
\begin{equation}\label{a_n}
|a_n|\leq \frac{1}{q_{n+1}} \leq \frac{1}{10 n^2q_n}.
\end{equation}
Hence, if we use the notation 
$$\overline{D}_{n,j}^1 = I_{n,j} \times [0,1] \subset D_{n,j}^1,$$
we have
\begin{eqnarray} \label{dbarre} R_{\a_{n+1}}^{m_n} (\overline{D}_{n,j}^1) \subset D_{n,j'}^2 \end{eqnarray}
for some $j' \in \ZZ$.

\subsubsection{Choice of the decompositions $\eta_n$.}

We define $\eta_n$ to be the partial decomposition  
of $M$ consisting of  the horizontal intervals $I_{n,j}\times\{r\} \subset D_{n,j}^1$, where
$r\in[1/(3n),1-1/(3n)]$, defined by (\ref{eqint_smooth}) and of the intervals $\overline{I}_{n,j} \times 
\{r\} $ with $r\in[1/(3n),1-1/(3n)]$ and 
$$\overline{I}_{n,j} = \left[\frac{j}{q_n}+ \frac{1}{2q_n} + \frac{1}{6nq_n} - a_n,
\frac{j+1}{q_n} -\frac{1}{6nq_n} - a_n \right].$$
It follows form (\ref{a_n}) that the intervals  $\overline{I}_{n,j} \times \{r\} $ are in $D_{n,j}^2$.

\begin{figure}
\psfrag{d}[][][0.7]{$D_{n,j}^1$}
\psfrag{dd}[][][0.7]{$D_{n,j}^2$}
\psfrag{d1}[][][0.7]{$D_{n,j}^2$}
\psfrag{dd1}[][][0.7]{$D_{n,j+1}^1$}
\psfrag{j}[][][0.7]{$I$}
\psfrag{fj}[][][0.7]{$\phi_n^{-1}(I)$}
\psfrag{j1}[][][0.7]{$I$}
\psfrag{fsj}[][][0.7]{$R_{\a_{n+1}}^{m_n} \phi_n^{-1}(I)$}
\psfrag{sj1}[][][0.7]{$R_{\a_{n+1}}^{m_n} \phi_n^{-1}(I)$}
\psfrag{fsfj}[][][0.7]{$\phi_n R_{\a_{n+1}}^{m_n} \phi_n^{-1}(I)$}
\psfrag{fsfj}[][][0.7]{$\ \Phi_n(I)$}
\psfrag{phi1}[][][0.7]{$\Phi_n(I)$}
\psfrag{ik}[][][0.7]{$\ \ I+[K]$}
\includegraphics[width=8cm]{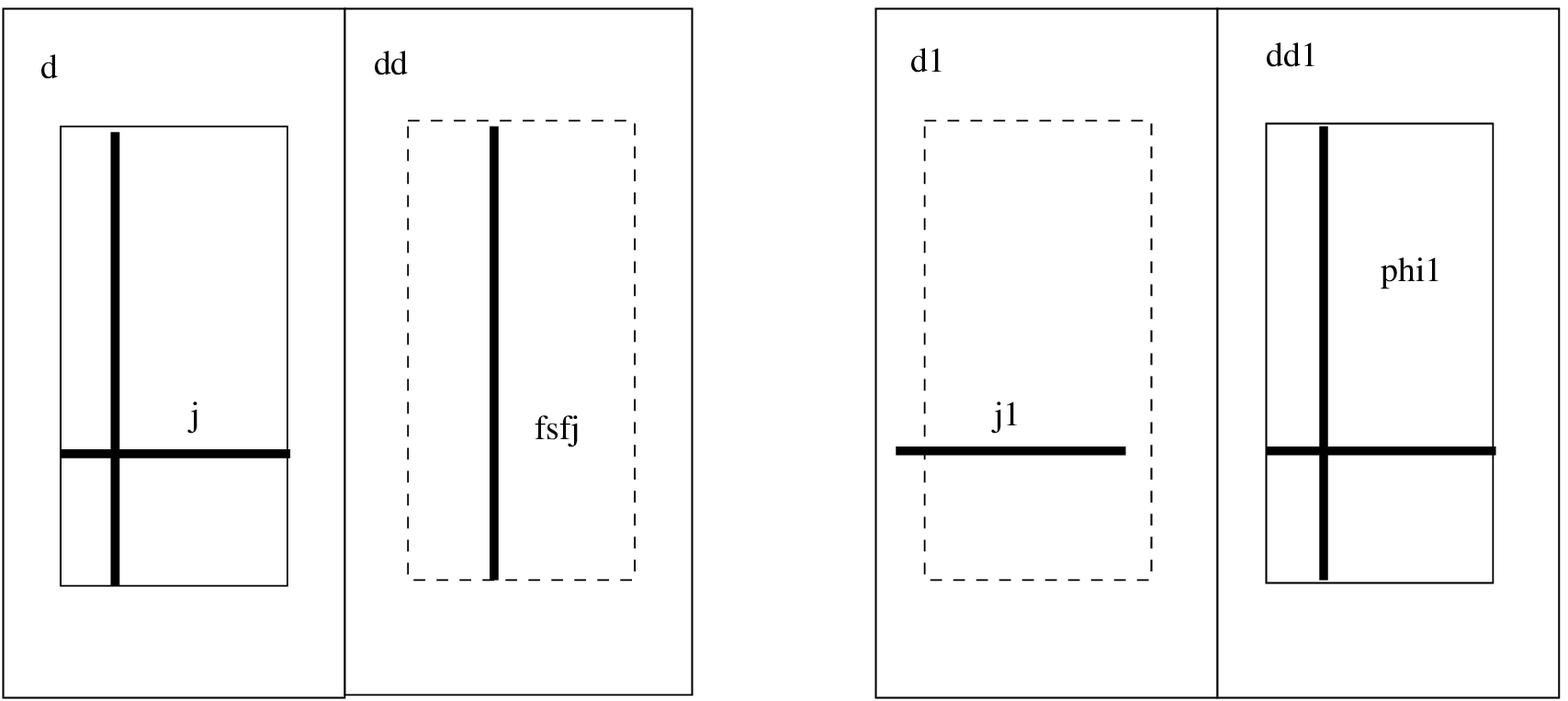}
\caption{Action of $\Phi_n$}\label{kartina3}
\end{figure}

\begin{lemm}\label{aff5}
The mapping  $\Phi_n= \phi_n \circ R_{\a_{n+1}}^{m_n}\circ  \phi_n^{-1}$ 
transforms the atoms of the
decomposition $\eta_n$ into  vertical intervals of the form 
$\{\theta\} \times [1/(3n),1-1/(3n)]$
for some~$\theta$. 
\end{lemm}
The proof is illustrated on Figure \ref{kartina3}.

\begin{proof}

Consider first an interval $I_n$ of the type $I_n=I_{n,j}\times \lbrace r \rbrace$,  $r\in 
[1/(3n),1-1/(3n)]$. By construction of $\phi_n$ (see \S \ref{ppphi}), we have that  $\phi_n^{-1}(I_{n})$ is 
a vertical segment of the form $\{\theta\} \times [1/(3n),1-1/(3n)]$
for some $\theta \in I_{n,j}$. From (\ref{dbarre}) we deduce that $R_{\a_{n+1}^{m_n}} \circ 
\phi_n^{-1}(I_{n}) = \{\theta'\} \times [1/(3n),1-1/(3n)] \subset  D_{n,j'}^2,$ for some $\theta' \in \TT$ 
and $j'\in \ZZ$ and we conclude using that $\phi_n$ acts as the identity on $D_{n,j'}^2.$

Similarly, for $r\in[1/(3n),1-1/(3n)]$ and an interval $I_n=\overline{I}_{n,j} \times \{r\} \in D_{n,j}^2$, 
we have that 
$$
\phi_n \circ R_{\a_{n+1}}^{m_n} \circ \phi_n^{-1}(I_n)=
\phi_n \circ R_{\a_{n+1}}^{m_n}(I_n)= \phi_n ({I}_{n,j'} \times \{r\}) = \{\theta\} 
\times[1/(3n),1-1/(3n)],$$
for some $j' \in \ZZ$ and $\theta \in \TT$.
\end{proof}


\subsubsection{Proof of Theorem \ref{thsmooth}.} Let the
diffeomorphisms $f_n$ be constructed as in  (\ref{eqsmooth_appr}),
following Lemma \ref{convergence} and  Lemma \ref{Nnda}, so that 
convergence of $f_n$, closeness to Identity of their limit $f$, as well as
(\ref{zvezdochka}) and (\ref{estimate}), hold. We want to apply
Proposition \ref{criterion} to get weak mixing. Since the sequence of
decompositions
$\eta_n\to \epsilon$ by construction, and since it consists of
intervals with length less than $1/q_n$, to finish it is enough to 
show that for any interval $I_n$ of the decomposition
$\eta_n$, and for 
$\Phi_n=\phi_n\circ R_{\a_{n+1}}^{m_n}\circ \phi_n^{-1}$, we have:
$\Phi_n(I_n)$ is  
$(0,2/(3n),0)$-distributed. 
The conditions of the definition follow immediately from the
construction and Lemma \ref{aff5}. Indeed, the projection of
$\Phi_n(I_n)$
to the $r$-axis is the interval $[1/(3n),1-1/(3n)]$, 
hence, in the definition of ($\ga,\de,\e$)-distribution 
(Definition \ref{aaaa}) we can take $\de=2/(3n)$. 
Furthermore, since the image of any interval 
$I_n$ is vertical, $\ga$ can be taken
equal to 0.
Finally, the restriction of ${\Phi_n}$ to ${I_n}$ being affine, one verifies that for any interval 
$\tilde{J}_n\subset J_n$:
$$
\la(I\cap\Phi_n^{-1}(\tilde J) )\la(J)=\la(I)\la(\tilde J).
$$
Hence, we take $\e=0$. 

We have verified the conditions
of Proposition \ref{criterion}. This implies weak mixing of the limit
diffeomorphism $f$.  \carre


\frenchspacing
\bibliographystyle{plain}

\end{document}